\documentclass[11pt, reqno, psamsfonts]{amsart}
\pdfoutput=1

\usepackage{amssymb}
\usepackage{amsthm}
\usepackage{amsmath}
\usepackage{latexsym}
\usepackage[T1]{fontenc}
\usepackage[utf8]{inputenc}
\usepackage[russian, french, english]{babel}

\usepackage{graphicx}
\usepackage{wrapfig}
\usepackage[justification=centering, labelfont=bf]{caption}
\usepackage{mathtools}
\usepackage{amsbsy}
\usepackage[inline]{enumitem}
\usepackage{mathrsfs}
\usepackage{array}
\usepackage{multicol}
\usepackage{stmaryrd}
\usepackage{cancel}
\usepackage{lmodern}
\usepackage{mathabx}
\usepackage{upgreek}
\usepackage{titlesec}
\usepackage{titletoc}
\usepackage[spacing=true,kerning=true,babel=true,tracking=true]{microtype}
\usepackage[shortcuts]{extdash}
\usepackage[foot]{amsaddr}
\usepackage[left=1in,right=1in,top=1in,bottom=1in,bindingoffset=0cm]{geometry}
\usepackage{bm}
\usepackage{centernot}
\usepackage{mdframed}
\usepackage[hidelinks]{hyperref}
\usepackage{xspace}
\usepackage{eucal,eufrak}
\usepackage[skins]{tcolorbox}
\usepackage{framed}
\usepackage[justification=centering, labelfont=bf]{caption}
\usepackage{tikz}
\usetikzlibrary{shapes,snakes}
\usetikzlibrary{arrows.meta}
\usetikzlibrary{decorations.pathmorphing}
\usetikzlibrary{patterns}
\usepackage{float}

\usepackage[
    backend=biber, style=alphabetic, sorting=nyt, maxnames=100,backref=true]{biblatex}
\title{\sffamily Borel versions of the Local Lemma and LOCAL algorithms for graphs of finite asymptotic separation index}
\date{}

\author{Anton~Bernshteyn}
\address{\normalfont (AB) Department of Mathematics, University of California, Los Angeles, USA}
\email{bernshteyn@math.ucla.edu}

\author{Felix~Weilacher}
\address{\normalfont (FW) Department of Mathematics, University of California, Berkeley, USA}
\email{weilacher@berkeley.edu}

\thanks{AB's research is partially supported by the NSF grant DMS-2045412 and the NSF CAREER grant DMS-2239187. FW's research is partially supported by the ARCS foundation, Pittsburgh chapter.}

\newtheoremstyle{bfnote}%
{}{}%
{\slshape}{}%
{\bfseries}{\bfseries.}%
{ }%
{\thmname{#1}\thmnumber{ #2}\thmnote{ \ep{\normalfont{}#3}}}

\theoremstyle{bfnote}
\newtheorem{theo}{Theorem}[section]
\newtheorem*{theo*}{Theorem}
\newtheorem{prop}[theo]{Proposition}
\newtheorem{lemma}[theo]{Lemma}
\newtheorem{claim}[theo]{Claim}

\newtheorem{corl}[theo]{Corollary}

\newtheorem*{corl*}{Corollary}

\theoremstyle{definition}
\newtheorem{defn}[theo]{Definition}
\newtheorem*{defn*}{Definition}
\newtheorem{exmp}[theo]{Example}
\newtheorem{exmps}[theo]{Examples}

\newtheorem{ques}[theo]{Question}
\newtheorem{remks}[theo]{Remarks}
\newtheorem*{exmp*}{Example}

\newtheorem{conj}[theo]{Conjecture}

\theoremstyle{remark}
\newtheorem*{ques*}{Question}
\newtheorem*{remk*}{Remark}

\newcommand*{\myproofname}{Proof}
\newenvironment{claimproof}[1][\myproofname]{\begin{proof}[#1]}{\end{proof}}

\makeatletter
\newcommand{\neutralize}[1]{\expandafter\let\csname c@#1\endcsname\count@}
\makeatother

\newenvironment{theocopy}[1]
{\renewcommand{\thetheo}{\ref{#1}}%
	\neutralize{theo}\phantomsection
	\begin{theo}}
	{\end{theo}}

\newcommand{\0}{\varnothing}
\newcommand{\set}[1]{\{#1\}}
\newcommand{\N}{{\mathbb{N}}}
\newcommand{\Z}{\mathbb{Z}}

\newcommand{\R}{\mathbb{R}}

\renewcommand{\P}{\mathbb{P}}
\newcommand{\E}{\mathbb{E}}
\renewcommand{\epsilon}{\varepsilon}

\renewcommand{\phi}{\varphi}
\renewcommand{\theta}{\vartheta}
\renewcommand{\leq}{\leqslant}
\renewcommand{\geq}{\geqslant}
\newcommand{\defeq}{\coloneqq}

\newcommand{\bemph}[1]{{\normalfont#1}} 
\newcommand{\ep}[1]{\bemph{(}#1\bemph{)}} 

\newcommand{\emphd}[1]{{\fontseries{b}\selectfont\textsf{#1}}}
\newcommand{\acts}{\mathrel{\reflectbox{$\righttoleftarrow$}}}

\newcommand{\G}{\Gamma}
\newcommand{\pto}{\dashrightarrow}
\newcommand{\dom}{\mathrm{dom}}

\newcommand{\rest}[2]{{{#1}\vert_{#2}}}
\newcommand{\B}{\mathscr{B}}
\newcommand{\pr}{\mathsf{p}}
\newcommand{\de}{\mathsf{d}}
\newcommand{\asi}{\mathsf{asi}}
\newcommand{\si}{\mathsf{si}}
\newcommand{\dist}{\mathsf{dist}}
\newcommand{\fs}[2]{{[{#1}]^{#2}}}
\newcommand{\finset}[1]{\fs{{#1}}{{<\infty}}}
\newcommand{\finfun}[2]{\fs{{#1} \pto {#2}}{{<\infty}}}
\newcommand{\LOCAL}{$\mathsf{LOCAL}$\xspace}
\newcommand{\Const}{\mathsf{Const}}
\newcommand{\pw}{polynomial\xspace}

\numberwithin{equation}{section}

\newenvironment{scproof}[1][]{\begin{proof}[\textsc{\upshape{Proof}}#1]}{\end{proof}}

\titleformat{\section}[block]{\large\bfseries\sffamily}{\thesection.}{1ex}{}
\titleformat{\subsection}[block]{\bfseries\sffamily}{\thesubsection.}{1ex}{}
\titleformat{\subsubsection}[block]{\itshape}{\bfseries\upshape\sffamily\thesubsubsection.}{1ex}{}

\titlespacing*{\section}{0pt}{*3}{*1}
\titlespacing*{\subsection}{0pt}{*3}{*1}
\titlespacing*{\subsubsection}{0pt}{*2}{*1}

\titlecontents{section}
 [1.5em] %
 {\smallskip}
 {\bfseries\thecontentslabel\hspace{1.02em}}
{\bfseries}
 {\,\,\titlerule*[0.77pc]{}\bfseries\contentspage}
\titlecontents{subsection}
 [4em] %
 {\smallskip}
 {\thecontentslabel\hspace{1.02em}}
{\hspace*{2.32em}}
 {\,\,\titlerule*[0.77pc]{.}\contentspage}

\renewbibmacro{in:}{}

\renewbibmacro*{volume+number+eid}{%
	\printfield{volume}%
	\setunit*{\addnbspace}
	\printfield{number}%
	\setunit{\addcomma\space}%
	\printfield{eid}}

\DeclareFieldFormat[article]{volume}{\textbf{#1}\space}
\DeclareFieldFormat[article]{number}{\mkbibparens{#1}}

\DeclareFieldFormat{journaltitle}{#1,}
\DeclareFieldFormat[thesis]{title}{\mkbibemph{#1}\addperiod}
\DeclareFieldFormat[article, unpublished, thesis]{title}{\mkbibemph{#1},}
\DeclareFieldFormat[book]{title}{\mkbibemph{#1}\addperiod}
\DeclareFieldFormat[unpublished]{howpublished}{#1, }

\DeclareFieldFormat{pages}{#1}

\DeclareFieldFormat[article]{series}{Ser.~#1\addcomma}

\setlist{topsep=3pt,itemsep=3pt}

\begin{document}\renewcommand{\thefootnote}{\fnsymbol{footnote}}


    \maketitle
    
    
    \begin{abstract}
        Asymptotic separation index is a parameter that measures how easily a Borel graph can be approximated by its subgraphs with finite components. In contrast to the more classical notion of hyperfiniteness, asymptotic separation index is well-suited for combinatorial applications in the Borel setting. The main result of this paper is a Borel version of the Lov\'asz Local Lemma---a powerful general-purpose tool in probabilistic combinatorics---under a finite asymptotic separation index assumption. As a consequence, we show that locally checkable labeling problems that are solvable by efficient randomized distributed algorithms admit Borel solutions on bounded degree Borel graphs with finite asymptotic separation index. From this we derive a number of corollaries, for example a Borel version of Brooks's theorem for graphs with finite asymptotic separation index.
    \end{abstract}
    
    \section{Introduction and main results}

    \subsection{Borel combinatorics and the asymptotic separation index}\label{sec:intro_asi}

    This paper is a contribution to the area of \emphd{descriptive combinatorics}, which investigates classical combinatorial notions---such as colorings, matchings, etc.---from the standpoint of descriptive set theory. Before we proceed, let us review some basic notation. We use $\N$ to denote the set of all non-negative integers and identify each $q \in \N$ with the $q$-element set $q = \set{i \in \N \,:\, i < q}$. For a set $X$ and a natural number $q \in \N$, we let
    \[
        [X]^q \,\defeq\, \set{A \subseteq X \,:\, |A| = q} \qquad \text{and} \qquad \finset{X} \,\defeq\, \set{A \subseteq X \,:\, \text{$A$ is finite}}.
    \]
    All graphs in this paper are undirected and simple. In other words, a graph $G$ consists of a vertex set $V(G)$ and an edge set $E(G) \subseteq [V(G)]^2$. When there is no possibility of confusion, we use the standard graph-theoretic convention and write $uv$ instead of $\set{u,v}$ to indicate an edge joining $u$ and $v$. Given $R \in \N$, an \emphd{$R$-ball} around a vertex $v$ in a graph $G$, denoted by $B_G(v, R)$, is the set of all vertices reachable from $v$ by a path of at most $R$ edges. 
    
    Throughout the paper, proper coloring of graphs will be used as a motivating example of a combinatorial problem.

    \begin{defn}[Colorings and the chromatic number]\label{defn:coloring}
        Let $G$ be a graph. Given $q \in \N$, a \emphd{proper $q$-coloring} of $G$ is a function $f \colon V(G) \to q$ such that $f(u) \neq f(v)$ for every edge $uv \in E(G)$. The \emphd{chromatic number} $\chi(G)$ of $G$ is the smallest $q \in \N$ such that $G$ has a proper $q$-coloring; if no such $q \in \N$ exists, we set $\chi(G) \defeq \infty$.\footnote{Although there is a rich and interesting theory concerning colorings with infinitely many colors (see, e.g., \cites{Komjath}[\S4]{KechrisMarks}{Geschke}), we will only focus on the case of finitely many colors in this paper.}
    \end{defn}
    
    
    Descriptive combinatorics applies graph-theoretic ideas to the study of {Borel graphs}:

    \begin{defn}[Borel graphs]\label{defn:Borel_graph}
        A \emphd{Borel graph} is a graph $G$ whose vertex set $V(G)$ is a standard Borel space and whose adjacency relation $\set{(u,v) \in V(G)^2 \,:\, uv \in E(G)}$ is a Borel subset of $V(G)^2$. 
    \end{defn}

    The reader is referred to \cite{KechrisDST,AnushDST} for descriptive set theory background, such as the definition of a standard Borel space\footnote{The reader unfamiliar with this terminology may assume that the space is, say, the unit interval $[0,1]$. This typically results in no loss of generality, thanks to the Borel isomorphism theorem \cite[Theorem 15.6]{KechrisDST}.}, and to \cite{KechrisMarks,Pikh_survey} for an overview of Borel graph theory. A useful observation is that if $X$ is a standard Borel space, then the sets $[X]^q$ for $q\in\N$ and $\finset{X}$ also carry natural standard Borel structures (see \S\ref{subsec:const}), and $G$ is a Borel graph if and only if $V(G)$ is a standard Borel space and $E(G)$ is a Borel subset of $[V(G)]^2$. 
    
    Since every countable subset of a standard Borel space is Borel, countable graphs are \ep{trivially} Borel, and hence Definition \ref{defn:Borel_graph} is only really interesting when $V(G)$ is uncountable. On the other hand, a large part of descriptive combinatorics deals with the case when $G$ is \emphd{locally countable} (or even \emphd{locally finite}), i.e., when every vertex of $G$ has countably (resp.\ finitely) many neighbors. Indeed, all graphs considered in this paper will be locally countable. Note that $G$ is locally countable if and only if every connected component of $G$ is countable. As a result, locally countable graphs studied in descriptive combinatorics usually have uncountably many components \ep{this is in contrast to classical, finite combinatorics, where it is common to only focus on connected graphs}.
    
    The fundamental question of descriptive combinatorics may be formulated as follows:
    \begin{ques}\label{ques:general}
        When can a combinatorial construction \ep{for example, a proper $q$-coloring} be performed on a Borel graph $G$ in a ``constructive'' fashion?
    \end{ques}
    The word ``constructive'' here can mean different things depending on the context. The most basic interpretation is to require all functions, sets, etc.\ appearing in the construction to be Borel. By introducing additional structure on $V(G)$, this requirement may be relaxed. For instance, one may ask for a construction that is measurable with respect to a measure $\mu$ on $V(G)$ or Baire-measurable with respect to a compatible Polish topology $\tau$ on $V(G)$. As a concrete example, we can apply this perspective to proper colorings by modifying Definition~\ref{defn:coloring} as follows:

    \begin{defn}[Borel, measurable, Baire-measurable chromatic numbers]\label{defn:Bcoloring}
        Let $G$ be a Borel graph. A proper $q$-coloring $f \colon V(G) \to q$ of $G$ is \emphd{Borel} if $f^{-1}(i)$ is a Borel subset of $V(G)$ for all $0 \leq i < q$. The \emphd{Borel chromatic number} $\chi_\mathsf{B}(G)$ of $G$ is the smallest $q \in \N$ such that $G$ has a Borel proper $q$-coloring; 
        if no such $q$ exists, we set $\chi_\mathsf{B}(G) \defeq \infty$. 
        
        Similarly, given a probability Borel measure $\mu$ or a compatible Polish topology $\tau$ on $V(G)$, we say that a proper $q$-coloring $f \colon V(G) \to q$ of $G$ is \emphd{$\mu$-measurable} or \emphd{$\tau$-Baire-measurable} if $f^{-1}(i)$ is a $\mu$-measurable, resp.~$\tau$-Baire-measurable subset of $V(G)$ for all $0 \leq i < q$. The \emphd{$\mu$-measurable chromatic number} $\chi_\mu(G)$ and the \emphd{$\tau$-Baire-measurable chromatic number} $\chi_\tau(G)$ are defined as the smallest $q \in \N$ such that $G$ has a $\mu$-measurable, resp.~$\tau$-Baire-measurable 
        proper $q$-coloring if such $q$ exists, and $\infty$ otherwise. 
    \end{defn}


    Definition~\ref{defn:Bcoloring} is due to Kechris, Solecki, and Todorcevic, who initiated the systematic study of Borel graphs and their combinatorics in the seminal paper \cite{KST}. Among other things, they demonstrated that greedy algorithms may be implemented on locally finite Borel graphs ``in a Borel way,'' leading, for instance, to the following result:

    \begin{theo}[{Kechris--Solecki--Todorcevic \cite[Proposition 4.6]{KST}}]\label{theo:KST}
        If $G$ is a Borel graph of finite maximum degree $\Delta$, then $\chi_\mathsf{B}(G) \leq \Delta + 1$.
    \end{theo}

    Unfortunately, later work revealed that without any additional assumptions, it is often impossible to improve on the initial results obtained in the Borel setting using greedy algorithms. The following theorem of Marks serves as a remarkable illustration of this, showing that the bound in Theorem~\ref{theo:KST} is optimal even if $G$ has no cycles \ep{graphs with no cycles are called \emphd{forests}}:

    \begin{theo}[{Marks \cite{Marks}}]\label{theo:Marks}
        For every $\Delta \in \N$, there exists a Borel forest $G$ of maximum degree $\Delta$ such that $\chi_\mathsf{B}(G) = \Delta + 1$.
    \end{theo}

    By contrast, every forest $G$ satisfies $\chi(G) \leq 2$. Furthermore, Theorem~\ref{theo:Marks} shows that the following classical upper bounds on the chromatic number fail in the Borel setting:

    \begin{theo}[{Brooks \cites{Brooks}[Theorem~5.2.4]{Diestel}}]\label{theo:Brooks}
        If $G$ is a graph of finite maximum degree $\Delta \geq 3$ and without a $(\Delta + 1)$-clique, then $\chi(G) \leq \Delta$.
    \end{theo}

    \begin{theo}[{Johansson \cites{Joh_triangle}[\S13]{MolloyReed}}]\label{theo:Johansson}
        There exist constants $C$, $\Delta_0 > 0$ such that if $G$ is a triangle-free graph of finite maximum degree $\Delta \geq \Delta_0$, then $\chi(G) \leq C\Delta/\log\Delta$.
    \end{theo}

    Both Brooks's and Johansson's theorems have analogs in the measurable/Baire-measurable setting \cites{MeasurableBrooks}[Theorem 3.9]{BernshteynDistributed}, showing that Marks's Theorem~\ref{theo:Marks} is a purely Borel phenomenon.

    In view of the negative results such as Theorem~\ref{theo:Marks}, it is natural to look for additional assumptions on a Borel graph $G$ that improve the behavior of its Borel combinatorics. In this regard, a helpful observation is that Borel combinatorics essentially trivialize on \emphd{component-finite} graphs, i.e., graphs in which every connected component is finite. It is well understood in the field that a combinatorial problem admits a Borel solution on a component-finite Borel graph $G$ if and only if it admits any solution on $G$ at all \ep{see \cite[\S5.3]{Pikh_survey} and \S\ref{subsec:cf} for two ways of making this statement precise}. For example, $\chi_\mathsf{B}(G) = \chi(G)$ for every component-finite Borel graph $G$. Therefore, one may hope that Borel combinatorics on a locally countable Borel graph $G$ would become more approachable if $G$ could in some sense be ``approximated'' by its component-finite subgraphs. 

    One well-studied notion that attempts to capture this intuition is {hyperfiniteness}. 
    As usual, a \emphd{subgraph} of a graph $G$ is a graph $H$ such that $V(H) \subseteq V(G)$ and $E(H) \subseteq E(G)$. We write $H \subseteq G$ to indicate that $H$ is a subgraph of $G$. If $G$ and $H$ are Borel graphs and $V(H) \subseteq V(G)$ is a Borel subset equipped with the relative $\sigma$-algebra of Borel sets, we say that $H$ is a \emphd{Borel subgraph} of $G$. Given a set $U \subseteq V(G)$, we let $G[U]$ denote the subgraph of $G$ \emphd{induced} by $U$, i.e., the graph with vertex set $U$ and edge set $E(G[U]) \defeq E(G) \cap [U]^2$. 
    Note that when $U \subseteq V(G)$ is a Borel set, $G[U]$ is a Borel subgraph of $G$. 
    
    \begin{defn}[Hyperfinite graphs]
        A Borel graph $G$ is \emphd{hyperfinite} if there is an increasing sequence $G_0 \subseteq G_1 \subseteq \cdots \subseteq G$ of component-finite Borel subgraphs of $G$ whose union is $G$.
    \end{defn}

    The systematic study of hyperfinite graphs was initiated by Weiss \cite{Weiss} and Slaman and Steel \cite{SS}, with important foundational work done by Dougherty, Jackson, and Kechris \cite{DJK} and Jackson, Kechris, and Louveau \cite{jackson2002countable}, among others. For an overview of this topic, see \cite{KechrisCBER}. While hyperfiniteness is in general a very powerful notion, with many applications to measurable constructions (see, e.g., \cite{CM16, BowenKunSabok}), counterexamples from \cite{HypComb} show that it is largely unhelpful for obtaining Borel solutions to combinatorial problems. For example, Marks's Theorem~\ref{theo:Marks} remains true under the hyperfiniteness assumption:

    \begin{theo}[{Conley--Jackson--Marks--Seward--Tucker-Drob \cite{HypComb}}]\label{theo:Marks_hyp}
        For every $\Delta \in \N$, there exists a hyperfinite Borel forest $G$ of maximum degree $\Delta$ such that $\chi_\mathsf{B}(G) = \Delta + 1$.
    \end{theo}

    In particular, Theorems~\ref{theo:Brooks} and \ref{theo:Johansson} fail in the Borel setting for hyperfinite graphs.

    Nevertheless, in the breakthrough paper \cite{Dimension}, Conley, Jackson, Marks, Seward, and Tucker-Drob succeeded in isolating a different and more subtle notion of ``approximating'' a graph by its component-finite subgraphs that turns out to be extremely useful in combinatorial arguments. For a graph $G$ and an integer $R \in \N$, $G^R$ denotes the graph with $V(G^R) \defeq V(G)$ in which two distinct vertices $u$, $v$ are adjacent if and only if they are joined by a path of at most $R$ edges in $G$.

    \begin{defn}[Asymptotic separation index]
        Let $G$ be a Borel graph. The \emphd{separation index} of $G$, in symbols $\si(G)$, is the smallest $s \in \N$ for which there exists a partition $V(G) = U_0 \sqcup \ldots \sqcup U_s$ of $V(G)$ into $s+1$ Borel subsets such that the induced subgraphs $G[U_0]$, \ldots, $G[U_s]$ are component-finite 
        (if no such $s \in \N$ exists, we set $\si(G) \defeq \infty$). 
        The \emphd{asymptotic separation index} of $G$, in symbols $\asi(G)$, is the supremum of $\si(G^R)$ taken over all $R \in \N$. 
    \end{defn}

    There is a growing body of literature demonstrating that various combinatorial problems admit Borel solutions on Borel graphs with finite asymptotic separation index \cite{BWKonig,Dimension,ASIalgorithms,FelixFinDim}. For example, it turns out that the Borel chromatic number of a Borel graph is bounded above by a function of its ordinary chromatic number and asymptotic separation index:

    \begin{theo}[{Conley--Jackson--Marks--Seward--Tucker-Drob \cite[Theorem 8.1]{Dimension}}]\label{theo:chi_bound}
        If $G$ is a Borel graph with $\asi(G) < \infty$ and $\chi(G) < \infty$, then $\chi_\mathsf{B}(G) \leq (\asi(G) + 1) \chi(G) - \asi(G)$.
    \end{theo}

    Another example is a result of Qian and the second named author, who showed in \cite{ASIalgorithms} that a Borel graph $G$ of maximum degree $\Delta \in \N$ and with $\asi(G) < \infty$ has a Borel proper edge-coloring with $\Delta + 1 + \asi(G)$ colors. Moreover, if $G$ is bipartite, this bound can be improved to $\Delta + \asi(G)$ by a result of Bowen and the second named author \cite{BWKonig}. (We discuss edge-colorings in more detail in \S\ref{subsubsec:edge-coloring}.) 

    Of course, these results are only impressive assuming there is a rich family of Borel graphs with finite asymptotic separation index. Thankfully, this is indeed the case, as many natural classes of graphs have asymptotic separation index at most $1$. Here is a non-exhaustive list of examples (the reader is referred to the cited papers for details and omitted definitions):

    \begin{theo}[{Conley--Jackson--Marks--Seward--Tucker-Drob \cite{Dimension}}]\label{theo:asi1_group}
        Suppose that $G$ is a Schreier graph of a free Borel action of a finitely generated group $\G$ on a standard Borel space. Then $\asi(G) \leq 1$, provided $\G$ has at least one of the following properties:
        \begin{itemize}
            \item $\G$ is virtually nilpotent;
            \item $\G$ is polycyclic;
            \item $\G = \Z_2 \wr \Z$, the lamplighter group;
            \item $\G = BS(1,2)$, the Baumslag--Solitar group;
            \item $\G$ is solvable and linear over $\mathbb{Q}$.
        \end{itemize}
    \end{theo}


    \begin{theo}[{AB--Yu \cite{BernshteynYu}}]
        If $G$ is a Borel graph of polynomial growth, then $\asi(G) \leq 1$.
    \end{theo}

    Furthermore, \emph{every} locally finite Borel graph has asymptotic separation index at most $1$ outside of a meager set of vertices. To state this formally, we say that a set $U \subseteq V(G)$ is \emphd{$G$-invariant} if there are no edges between $U$ and $V(G)\setminus U$ (i.e., if $U$ is a union of connected components of $G$). As usual, we call a set in a topological space \emphd{comeager} if its complement is meager.

    \begin{theo}[{Conley--Jackson--Marks--Seward--Tucker-Drob \cite[Theorem 4.8(b)]{Dimension}}] \label{theo:ReduceAsiBM}
        Let $G$ be a locally finite Borel graph. Fix a compatible Polish topology $\tau$ on $V(G)$. Then there exists a $G$-invariant comeager Borel set $U \subseteq G$ such that $\asi(G[U]) \leq 1$.
    \end{theo}

    And if $G$ is hyperfinite, then similarly $\asi(G) \leq 1$ on a set of full measure:

    \begin{theo}[{\cite{FelixMeasure}}]\label{theo:ReduceAsiMeasure}
        Let $G$ be a hyperfinite locally finite Borel graph. Fix a probability Borel measure $\mu$ on $V(G)$. Then there is a $G$-invariant $\mu$-conull Borel set $U \subseteq G$ such that $\asi(G[U]) \leq 1$.
    \end{theo}

    Both Theorems \ref{theo:ReduceAsiBM} and \ref{theo:ReduceAsiMeasure} are essentially present in the paper \cite{CM16} by Conley and Miller, albeit phrased without the term ``asymptotic separation index.''
    
    It is not a coincidence that in all the above examples, the asymptotic separation index is not only finite but in fact bounded by $1$: it is not known whether there exist any Borel graphs $G$ with $1 < \asi(G) < \infty$ \cite[15]{Dimension}. The problem of whether such graphs exist is especially interesting in view of results such as Theorem~\ref{theo:chi_bound} as well as some of our main results below.

    The aim of this paper is to  make several additions to the list of combinatorial constructions that can be performed in a Borel way on graphs with finite asymptotic separation index. For example, we establish a Borel version of Brooks's Theorem~\ref{theo:Brooks} for Borel graphs $G$ with $\asi(G) < \infty$:

    \begin{tcolorbox}
    \begin{theo}[Borel Brooks for finite $\asi$]\label{theo:BorelBrooks}
        Let $G$ be a Borel graph with finite maximum degree $\Delta \geq 3$ and without a $(\Delta + 1)$-clique. If $\asi(G) < \infty$, then $\chi_\mathsf{B}(G) \leq \Delta$.
    \end{theo}
    \end{tcolorbox}
    
    We also establish a Borel version of Johansson's Theorem~\ref{theo:Johansson} for Borel graphs $G$ with $\asi(G) < \infty$. A result of this type can already be derived from Theorem~\ref{theo:chi_bound}. Indeed, by combining the best currently known bound on the chromatic number of triangle-free graphs (due to Molloy \cite{Molloy}) with Theorem~\ref{theo:chi_bound}, we get the following result for triangle-free Borel graphs $G$ with $\asi(G) < \infty$:
    \[
        \chi_\mathsf{B}(G) \,\leq\, (1 + \epsilon)(\asi(G) + 1) \frac{\Delta}{\log\Delta} \,-\, \asi(G),
    \]
    for any given $\epsilon > 0$ and large enough $\Delta$. In other words, this approach yields a factor in front of $\Delta/\log \Delta$ that is close to $\asi(G) + 1$. We prove a bound with a factor in front of $\Delta/\log \Delta$ that is
    independent of $\asi(G)$. Furthermore, when $G$ has no cycles of length $3$ and $4$, our bound matches the best known bound on the ordinary chromatic number of $G$ (due to Kim \cite{Kim95}):
    
    \begin{tcolorbox}
    \begin{theo}[Borel colorings of graphs without short cycles and with finite $\asi$]\label{theo:BorelJohansson}
            For every $\epsilon > 0$, there is $\Delta_0 > 0$ with the following property. Let $G$ be a Borel graph of finite maximum degree $\Delta \geq \Delta_0$ with $\asi(G) < \infty$. If $G$ is triangle-free, then $\chi_\mathsf{B}(G) \leq (4 + \epsilon)\Delta/\log\Delta$, and if $G$ has no $3$- and $4$-cycles, then $\chi_\mathsf{B}(G) \leq (1 + \epsilon)\Delta/\log\Delta$.
    \end{theo}
    \end{tcolorbox}
    
    It remains an interesting open question whether the $4$ in the above theorem can be replaced by $1$ in the triangle-free case.
    
    While most work in descriptive combinatorics addresses specific combinatorial problems, several recent papers attack the fundamental Question~\ref{ques:general} for relatively general classes of problems \cite{BerShift, Ber_Baire, BernshteynDistributed, Ber_cont, trees, CGMPT, paths, GRgrids, ASIalgorithms, FelixComput}. We continue this line of research here by establishing a Borel version of a powerful and versatile tool from probabilistic combinatorics, the so-called {Lov\'asz Local Lemma}, for Borel graphs with finite asymptotic separation index. Statements such as Theorems~\ref{theo:BorelBrooks} and \ref{theo:BorelJohansson} are among the consequences of our main result, as we explain in \S\ref{subsec:intro_appl}. 

    \subsection{The Lov\'asz Local Lemma and its Borel versions}\label{subsec:intro_Borel_LLL} 

     Typical combinatorial problems require assigning a ``color'' \ep{an element of some finite set} to every member of a given structure \ep{e.g., to every vertex of a graph} in a way that fulfills a prescribed set of constraints. This idea is formally captured in the following definition:
    
    \begin{defn}[{Constraint satisfaction problems}]\label{defn:CSP}
		Let $q$ be a positive integer. By a \emphd{$q$-coloring} of a set $S$ we mean a function $f \colon S \to q$.
			
		Given a set $X$ and a finite subset $D \subseteq X$, an \emphd{$(X,q)$-constraint} \ep{or simply a \emphd{constraint} if $X$ and $q$ are understood} with \emphd{domain} $D$ is a set $B \subseteq q^D$ of $q$-colorings of $D$. If $B$ is a constraint with domain $D$, then we write $\dom(B) \defeq D$. 
        We say that a $q$-coloring $f \colon X \to q$ of $X$ \emphd{violates} a constraint $B$ 
        if $\rest{f}{\dom(B)} \in B$, and \emphd{satisfies} $B$ otherwise.
			
		A \emphd{constraint satisfaction problem} \ep{a \emphd{CSP} for short} $\B$ \emphd{on} a set $X$ with \emphd{range} $q$, in symbols \[\B \colon X \to^? q,\] is a set of $(X,q)$-constraints. A \emphd{solution} to a CSP $\B \colon X \to^? q$ is a $q$-coloring $f \colon X \to q$ that satisfies every constraint $B \in \B$. We say that $\B$ is \emphd{satisfiable} if it has a solution.
	\end{defn}
	
	In other words, each constraint in a CSP $\B \colon X \to^? q$ is interpreted as a set of finite ``forbidden'' (or ``bad'') patterns that are not allowed to appear in a solution $f \colon X \to q$.
 
    \begin{exmp}[Proper coloring as a CSP]\label{exmp:coloring}
        Finding a proper $q$-coloring of a graph $G$ can be represented by a CSP $\B \colon V(G) \to^? q$ as follows. For an edge $uv \in E(G)$, let $B_{uv}$ be the constraint with domain $\set{u,v}$ given by $B_{uv} \defeq \set{(u \mapsto i,\, v \mapsto i) \,:\, 0 \leq i < q}$. Then a coloring $f \colon V(G) \to q$ satisfies $B_{uv}$ precisely when $f(u) \neq f(v)$. Therefore, $f$ is a solution to $\B \defeq \set{B_{uv} \,:\, uv \in E(G)}$ if and only if it is a proper $q$-coloring of $G$.
    \end{exmp}

    The \emphd{Lov\'asz Local Lemma} (the \emphd{LLL} for short) provides a sufficient condition guaranteeing that a CSP $\B \colon X \to^? q$ is satisfiable. To state it, we associate to $\B$ two numerical parameters, $\pr(\B)$ and $\de(\B)$, defined as follows. For a constraint $B \in \B$, the \emphd{probability} of $B$ is the quantity
	\[
	    \P[B] \,\defeq\, \frac{|B|}{q^{|\dom(B)|}}.
	\]
	In other words, $\P[B]$ is the probability that $B$ is violated by a uniformly random $q$-coloring. Let
	\[
	    \pr(\B) \,\defeq\, \sup_{B \in \B} \P[B].
	\]
	The \emphd{maximum dependency degree} $\de(\B)$ of $\B$ is the value
	\[
	    \de(\B) \,\defeq\, \sup_{B \in \B} |\set{B' \in \B \,:\, B \neq B' \text{ and } \dom(B') \cap \dom(B) \neq \0}|.
	\]

	\begin{lemma}[{Lov\'asz Local Lemma \cites{EL}{SpencerRamsey}[Corollary 5.1.2]{AS}}]\label{lemma:LLL}
		If $\B$ is a CSP such that
		\[
			\pr(\B) \, (\de(\B) + 1) \,\leq\, 1/e,
		\]
		where $e = 2.71\ldots$ is the base of the natural logarithm, then $\B$ has a solution.
	\end{lemma}
    
    \begin{remks}
    \begin{enumerate}[wide, label=\ep{\normalfont\roman*}]
        \item The statement of the LLL given above is somewhat less general than the more classical version usually found in the combinatorics literature. Specifically, the probabilities of the constraints in Lemma \ref{lemma:LLL} are computed by considering a \emph{uniformly random} $q$-coloring $f \colon X \to q$; that is, we view $(f(x) \,:\, x \in X)$ as a collection of mutually independent discrete random variables, each distributed uniformly over a $q$-element set. This specialized set-up for the LLL is called the \emphd{variable version} of the LLL (the name is due to Kolipaka and Szegedy \cite{KolipakaSzegedy}). Even though this setting is not the most general, it does encompass virtually all standard applications and is often viewed as the ``right one'' for algorithmic considerations (see, e.g., \cite{Beck,MT,FG}). For the statement of the LLL in abstract probability spaces, see \cite[\S5.1]{AS}.

        \item The LLL is usually stated and proved in the case when the CSP $\B$ comprises only finitely many constraints. 
        Nevertheless, a routine compactness argument shows that an infinite CSP $\B$ is satisfiable if and only if all its finite subsets $\B' \subset \B$ are satisfiable \ep{see, e.g., \cite[proof of Theorem 5.2.2]{AS}}, so the conclusion of Lemma~\ref{lemma:LLL} for infinite $\B$ follows from the finite case.
    \end{enumerate}
    \end{remks}

    The LLL is used throughout combinatorics; for a sample of its applications, see \cites[\S5]{AS}{MolloyReed}.
    
    
    Here we are interested in the LLL from the perspective of descriptive combinatorics. This line of research forms part of a broader program concerning versions of the LLL that are ``constructive'' in various senses: algorithmic \cite{Beck,BGRDeterministicLLL,FG,MT}, computable \cite{RSh}, or---as in this paper---Borel/measurable \cite{BerShift,BernshteynDistributed, Ber_cont, CGMPT}. 
    
    Let $X$ be a standard Borel space and let $q \in \N^+$. Then the set $\Const(X,q)$ of all $(X,q)$-constraints also carries a natural standard Borel structure (see \S\ref{subsec:const}), so we can define Borel CSPs as follows:

    \begin{defn}[Borel CSPs]\label{defn:BorelCSP}
        Let $X$ be a standard Borel space and let $q \in \N^+$. A CSP $\B \colon X \to^? q$ is \emphd{Borel} if it is a Borel subset of $\Const(X,q)$.
    \end{defn}

    Heuristically, a CSP $\B \colon X \to^? q$ is Borel as long as the constraints that need to be satisfied are specified ``explicitly.''

    \begin{exmp}[Proper coloring as a Borel CSP]\label{exmp:coloringBorelCSP}
        If $G$ is a Borel graph, then the CSP from Example~\ref{exmp:coloring} that encodes proper $q$-coloring problem on $G$ is easily seen to be Borel.
    \end{exmp}
    
    
    
    We now ask the following natural question:
    
    \begin{ques}\label{ques:definableLLL}
        Given a Borel CSP $\B \colon X \to^? q$ on a standard Borel space $X$, what LLL-style conditions guarantee that $\B$ has a Borel solution $f \colon X \to q$?
    \end{ques}
    
    Without any additional assumptions on $\B$, Question~\ref{ques:definableLLL} has recently been fully resolved in a series of contributions by several authors. In the following statement and in the remainder of the paper, we say that a CSP $\B$ is \emphd{bounded} if $\sup_{B \in \B} |\dom(B)| < \infty$.\footnote{The boundedness assumption can likely be removed, but it makes the results of \cite{BernshteynDistributed} easier to apply and is satisfied in most applications.}
    
    \begin{theo}[Brandt--Grunau--Rozho{\v{n}} \cite{BGRDeterministicLLL}, AB \cite{BernshteynDistributed}, Thornton \cite{ThorntonOrientation}]\label{theo:BorelLLL}
        \mbox{}
        
        \begin{enumerate}[label=\ep{\upshape{\roman*}}]
            \item\label{item:positive} Let $\B \colon X \to^? q$ be a bounded Borel CSP on a standard Borel space $X$. If
            \[
                \pr(\B) \, 2^{\de(\B)} \, < \, 1,
            \]
            then $\B$ has a Borel solution $f \colon X \to q$.
            
            \item\label{item:negative} On the other hand, for every $d \in \N$, there exists a bounded Borel CSP $\B$ on a standard Borel space such that $\de(\B) = d$ and $\pr(\B) = 2^{-d}$, so
            \[
                \pr(\B) \, 2^{\de(\B)} \, = \, 1,
            \]
            yet $\B$ has no Borel solution.
        \end{enumerate}
    \end{theo}

    Let us say a few words about how Theorem~\ref{theo:BorelLLL} follows from the cited sources. In \cite{BGRDeterministicLLL}, Brandt, Grunau, and Rozho\v{n} designed an efficient deterministic distributed algorithm for finding solutions to CSPs under the condition $\pr\, 2^{\de} < 1$ \ep{their work builds on earlier contributions of Brandt, Maus, and Uitto \cite{BMUDeterministicLLL}}. In \cite{BernshteynDistributed}, the first named author established certain general results that allow using efficient distributed algorithms to derive conclusions in descriptive combinatorics. In particular Theorem~\ref{theo:BorelLLL}\ref{item:positive} follows immediately from the Brandt--Grunau--Rozho\v{n} algorithm via \cite[Theorem 2.10]{BernshteynDistributed}. 
    \ep{See \S\ref{subsec:dist} for a more detailed discussion of the connections between distributed computing and descriptive combinatorics.}
    
    The second part of Theorem~\ref{theo:BorelLLL} is proved by analyzing the so-called \emphd{sinkless orientation problem} on regular graphs. Let $d \in \N$ and let $G$ be a $d$-regular graph. An orientation of $G$ is \emphd{sinkless} if the outdegree of every vertex is positive. A sinkless orientation of $G$ can be naturally encoded as a solution to a certain CSP $\B_{\text{sinkless}} \colon E(G) \to^? 2$. Here the color of each edge $e \in E(G)$ indicates the direction in which $e$ is oriented. For example, we may fix a function $c \colon E(G) \to V(G)$ sending each edge to one of its endpoints and say that an edge $e$ is oriented towards $c(e)$ if and only if its color is $0$. Then we define, for each vertex $v$, a constraint $B_v$ with domain $\dom(B_v) \defeq \set{e \in E(G) \,:\, e \ni v}$ that requires the outdegree of $v$ to be positive, and let \[\B_{\text{sinkless}} \,\defeq\, \set{B_v \,:\, v \in V(G)}.\] Clearly, we have $\de(\B_{\text{sinkless}}) = d$ and $\pr(\B_{\text{sinkless}}) = 2^{-d}$. Furthermore, if $G$ is a Borel graph, then the CSP $\B_{\text{sinkless}}$ can be made Borel by choosing a Borel function $c \colon E(G) \to V(G)$ using the Luzin--Novikov theorem \cite[Theorem 18.10]{KechrisDST}. 
    Nevertheless, Thornton \cite[Theorem 3.5]{ThorntonOrientation} applied the determinacy method of Marks \cite{Marks} to construct, for each $d \in \N$, a Borel $d$-regular forest $G$ that does not admit a Borel sinkless orientation. This verifies Theorem~\ref{theo:BorelLLL}\ref{item:negative}.
    
    While Theorem~\ref{theo:BorelLLL} provides a complete characterization of the range of values for the parameters $\pr$ and $\de$ under which a Borel version of the LLL holds, its downside is that the inequality $\pr\, 2^\de < 1$ is much harder to satisfy than 
    the usual LLL condition $\pr (\de + 1) \leq 1/e$. In fact, 
    we are not aware of any applications of substantial combinatorial interest where the bound $\pr\, 2^\de < 1$ holds and Theorem~\ref{theo:BorelLLL}\ref{item:positive} can be invoked. This motivates the search for other Borel versions of the LLL which 
    make extra assumptions on the structure of the CSP in question. 
    To formulate such assumptions, it is convenient to associate to each CSP a graph:

    \begin{defn}[Graphs associated to CSPs]
        Let $\B \colon X \to^? q$ be a CSP. Then $G_\B$ is the graph with vertex set $X$ in which two vertices $x$, $y \in X$ are adjacent if and only if $x \neq y$ and there is some constraint $B \in \B$ with $x$, $y \in \dom(B)$.
    \end{defn}
    
    Note that if $\B \colon X \to^? q$ is a Borel CSP and $\de(\B)$ is finite, 
    then $G_\B$ is a locally finite Borel graph. Moreover, if additionally $\B$ is bounded, then the maximum degree of $G_\B$ is finite.

    An important result of Cs{\'o}ka, Grabowski, M{\'a}th{\'e}, Pikhurko, and Tyros is a Borel version of the LLL for CSPs whose associated graphs have subexponential growth. Here we say that a graph $G$ is of \emphd{subexponential growth} if for every $\epsilon > 0$ there is $r > 0$ such that for all $R \geq r$ and all $v \in V(G)$, the $R$-ball around $v$ in $G$ contains at most $e^{\epsilon R}$ vertices.

    \begin{theo}[{Borel LLL for subexponential growth; Cs{\'o}ka--Grabowski--M{\'a}th{\'e}--Pikhurko--Tyros \cite[Theorem 4.5]{CGMPT}}]\label{theo:subexp}
        Let $\B \colon X \to^? q$ be a Borel CSP on a standard Borel space $X$. Suppose that the graph $G_\B$ is of subexponential growth. If 
        \[ \pr(\B)(\de(\B) + 1) \,\leq\, 1/e,\] then $\B$ has a Borel solution $f \colon X \to q$.
    \end{theo}

    
    We investigate the situation where instead of subexponential growth, the graph $G_\B$ has finite asymptotic separation index. Actually, for our first main result, we only need an upper bound on the \emph{separation index} of $G_\B$ rather than its \emph{asymptotic} separation index. (Even more precisely, we use a bound on a certain auxiliary parameter closely related to $\si(G_\B)$, which we call the shattering number of $\B$; see \S\ref{subsec:proof_main_easy} for details.) 
    
    \begin{tcolorbox}
    \begin{theo}[Borel LLL with bounded separation index]\label{theo:main_easy_asi}
        Let $\B \colon X \to^? q$ be a Borel CSP on a standard Borel space $X$. Suppose that $\si(G_\B) \leq s < \infty$. If
        \[
            \pr(\B) \, (\de(\B) + 1)^{s+1} \,\leq\, e^{-s-1},
        \]
        then $\B$ has a Borel solution $f \colon X \to q$. 
    \end{theo}
    \end{tcolorbox}

     It is a standard observation (see, e.g., \cite{Beck, BernshteynDistributed, CP, FG, GHK}) that 
     for a majority of applications, instead of the usual LLL condition $\pr(\de+1)\leq 1/e$, it is enough to have a version of the LLL that holds under a \emphd{\pw criterion}, i.e., with a bound of the form $\pr \, f(\de) \leq 1$ for some polynomial $f$. When $s$ is treated as a constant parameter, our Theorem~\ref{theo:main_easy_asi} is precisely of this type. (This should be contrasted with the exponential bound $\pr \, 2^\de < 1$ in Theorem~\ref{theo:BorelLLL}.) 
    For instance, if $G_\B$ has asymptotic separation index at most $1$ (which happens in the numerous situations discussed in \S\ref{sec:intro_asi}), we can apply the trivial inequality $\si(G_\B) \leq \asi(G_\B)$ to obtain the following special case of Theorem~\ref{theo:main_easy_asi}: 

    \begin{tcolorbox}
    \begin{corl}[Borel LLL for $\asi \leq 1$]\label{corl:main_asi1}
        Let $\B \colon X \to^? q$ be a Borel CSP on a standard Borel space $X$. Suppose that $\asi(G_\B) \leq 1$. If
        \[
            \pr(\B) \, (\de(\B) + 1)^2 \,\leq\, e^{-2},
        \]
        then $\B$ has a Borel solution $f \colon X \to q$. 
    \end{corl}
    \end{tcolorbox}

    Most of the examples listed in Theorem~\ref{theo:asi1_group} have or can have exponential growth rate. In such cases, Corollary~\ref{corl:main_asi1} is an advancement over both Theorems~\ref{theo:subexp} and \ref{theo:BorelLLL}\ref{item:positive}. Interestingly, it is not known whether there exist any Borel graphs $G$ of subexponential growth with $\asi(G) > 1$. 
    
    Of course, the result of Corollary~\ref{corl:main_asi1} holds for other finite values of $\asi(G_\B)$ as well, but with the exponent $2$ replaced by $\asi(G_\B) + 1$. Our second main result is a Borel version of the LLL for bounded CSPs $\B$ satisfying $\asi(G_\B) < \infty$, which holds under a \pw criterion independent of the actual value of $\asi(G_\B)$:

    \begin{tcolorbox}
    \begin{theo}[Borel LLL for finite $\asi$]\label{theo:main_speedup}
        Let $\B \colon X \to^? q$ be a bounded Borel CSP on a standard Borel space $X$. Suppose that $\asi(G_\B) < \infty$. If
        \[
            \pr(\B) \, (\de(\B) + 1)^8 \,\leq\, 2^{-15},
        \]
        then $\B$ has a Borel solution $f \colon X \to q$. 
    \end{theo}
    \end{tcolorbox}

    The above results add further impetus to solve the problem of whether there exist Borel graphs $G$ with $1 < \asi(G) < \infty$, since if they do not, Corollary~\ref{corl:main_asi1} would supersede Theorem~\ref{theo:main_speedup}. 


    In addition to Borel versions of the LLL, we can also consider its measurable and Baire-measurable versions. To date, the most general result along these lines is the following:

    \begin{theo}[{AB \cite[Theorem 2.20]{BernshteynDistributed}}]\label{theo:meas_LLL_old}
        Let $\B \colon X \to^? q$ be a bounded Borel CSP on a standard Borel space $X$. Suppose that
		\[
            \pr(\B) \, (\de(\B) + 1)^8 \,\leq\, 2^{-15}.
        \]
		Then the following conclusions hold:
		\begin{enumerate}[label=\ep{\normalfont\roman*}]
			\item\label{item:LLLmeas} If $\mu$ is a probability Borel measure on $X$, then $\B$ has a $\mu$-measurable solution.
			
			\item\label{item:LLLBaire} If $\tau$ is a compatible Polish topology on $X$, then $\B$ has a $\tau$-Baire-measurable solution.
		\end{enumerate}
    \end{theo}

    Thus, the LLL holds measurably/Baire-measurably with a \pw criterion. However, it is still not known whether the usual criterion $\pr (\de + 1) \leq 1/e$ suffices (we conjecture that it does). By combining Corollary~\ref{corl:main_asi1} with Theorem~\ref{theo:ReduceAsiBM}, we immediately obtain a quantitative improvement to Theorem~\ref{theo:meas_LLL_old} in the Baire-measurable case:

    \begin{tcolorbox}        
    \begin{corl}[Baire-measurable LLL]\label{corl:BMLLL_asi}
        Let $\B \colon X \to^? q$ be a Borel CSP on a standard Borel space $X$. Suppose that
		\[
            \pr(\B) \, (\de(\B) + 1)^2 \,\leq\, e^{-2}.
        \]
        If $\tau$ is a compatible Polish topology on $X$, then $\B$ has a $\tau$-Baire-measurable solution. 
    \end{corl}
    \end{tcolorbox}

    Similarly, if $G_\B$ is hyperfinite, Theorem~\ref{theo:ReduceAsiMeasure} yields an improvement in the measurable setting:

    \begin{tcolorbox}        
    \begin{corl}[Measurable hyperfinite LLL]\label{corl:MLLL_asi}
        Let $\B \colon X \to^? q$ be a Borel CSP on a standard Borel space $X$. Suppose that the graph $G_\B$ is hyperfinite and that
		\[
            \pr(\B) \, (\de(\B) + 1)^2 \,\leq\, e^{-2}.
        \]
        If $\mu$ is a probability Borel measure on $X$, then $\B$ has a $\mu$-measurable solution.
    \end{corl}
    \end{tcolorbox}

    \subsection{Distributed algorithms}\label{subsec:dist}

    A recent trend in descriptive combinatorics, initiated by the first named author in \cite{BernshteynDistributed}, is to investigate the surprisingly intimate relationship between descriptive combinatorics and \emphd{distributed computing}---an area of computer science concerned with problems that can be solved efficiently by a decentralized network of processors. It turns out that the existence of an efficient distributed algorithm that solves a combinatorial problem on a class of finite graphs often implies that the same problem admits Borel, measurable, etc.\ solutions on certain Borel graphs \cite{BernshteynDistributed}, and sometimes the reverse implication also holds \cite{Ber_cont,trees,GRgrids}. Using Theorem~\ref{theo:main_easy_asi}, we contribute a new result to this avenue of research, which allows using randomized distributed algorithms to solve combinatorial problems in a Borel way on graphs with finite asymptotic separation index.

    The class of combinatorial problems we consider here are the so-called locally checkable labeling \ep{LCL} problems, which were isolated by Naor and Stockmeyer in \cite{NaorStock} (in \cite{BernshteynDistributed} the term ``local coloring problems'' is used instead). The definition of an LCL problem is very similar to the definition of a CSP given in Definition~\ref{defn:CSP}; the main difference is that an LCL problem is defined on an underlying graph using the graph's structure:

    \begin{defn}[Locally checkable labeling problems]\label{defn:LCL}
        A \emphd{labeled graph} is a pair $(G, f)$, where $G$ is a graph and $f \colon V(G) \to \N$ is a function, called a \emphd{labeling}. A \emphd{rooted labeled graph} is a triple $(G, f, v)$, where $(G,f)$ is a labeled graph and $v \in V(G)$ is a vertex, called a \emphd{root}. Two rooted labeled graphs $(G,f,v)$ and $(G',f',v')$ are \emphd{isomorphic} if there is a graph isomorphism between $G$ and $G'$ that sends $v$ to $v'$ and preserves the labeling. Given a labeled graph $(G,f)$, a vertex $v \in V(G)$, and an integer $R \in \N$, we let $[G, f, v]_R$ be the isomorphism type of the rooted labeled graph
        \[
            \big(G[B_G(v,R)], \,\rest{f}{B_G(v,R)}, \, v\big),
        \]
        where $B_G(v,R)$ is the $R$-ball around $v$ in $G$.
        
        A \emphd{locally checkable labeling problem} (an \emphd{LCL problem} for short) of \emphd{radius} $R$ is a pair $\Pi = (R, \mathcal{P})$, where $R \in \N$ and $\mathcal{P}$ is a map that sends isomorphism classes of finite rooted labeled graphs to $\set{0, 1}$. If $\Pi = (R, \mathcal{P})$ is an LCL problem and $(G,f)$ is a locally finite labeled graph, then $f$ \emphd{solves} $\Pi$ on $G$, or is a \emphd{$\Pi$-coloring} of $G$, if $\mathcal{P}([G, f, v]_R) = 1$ for all $v \in V(G)$. A graph $G$ is \emphd{$\Pi$-colorable} if it has a $\Pi$-coloring $f \colon V(G) \to \N$.
    \end{defn}

    Informally, an LCL problem $\Pi = (R, \mathcal{P})$ is a ``rule'' that decides whether a labeling $f$ of $G$ is ``valid'' by looking at the restrictions of $f$ to $R$-balls around individual vertices.

    \begin{exmp}[{Proper coloring as an LCL problem}]\label{exmp:coloringLCL}
		A typical example of an LCL problem is proper $q$-coloring, since whether a coloring of a graph $G$ is proper is determined by its restrictions to $1$-balls in $G$. Explicitly, given \ep{the isomorphism type of} a finite rooted labeled graph $(G,f,v)$, set $\mathcal{P}(G, f, v) \defeq 1$ if and only if $f(v) < q$ 
        and $f(u) \neq f(v)$ for all neighbors $u$ of $v$. If we let $\Pi \defeq (1, \mathcal{P})$, then a $\Pi$-coloring of a locally finite graph $G$ is the same as a proper $q$-coloring of $G$.
	\end{exmp}

    The model of distributed computation relevant to our work is called \LOCAL\footnote{This is not an acronym.}. It was introduced by  Linial in \cite{Linial} (although there are some earlier related results, e.g., by Alon, Babai, and Itai \cite{ABI}, Luby \cite{Luby}, and Goldberg, Plotkin, and Shannon \cite{GPSh}). For an introduction to this subject, see the book \cite{BE} by Barenboim and Elkin.

    The \LOCAL model is intended to quantify the difficulty of transforming local data into a global solution to a problem on a large (but finite!) decentralized communication network. Informally, we imagine each vertex of an $n$-vertex graph $G$ to be occupied by a processor that may pass messages to its neighbors. The length of the messages is unrestricted, and the computational power available to each processor is unlimited (in other words, any computations carried out by individual processors take a single unit of time). Eventually, every processor must decide on its own part of the output; in the context of LCL problems, each vertex must decide on its own label. Since the processors may only communicate along the edges of $G$, in $T$ time units, a processor can only collect information from its $T$-ball in $G$---this is why the model is called ``\LOCAL''. In effect, a \LOCAL algorithm with running time $T$ can be thought of as a function that determines the output at each vertex of $G$ based only on the isomorphism type of the $T$-ball around it \cite[\S4.1.2]{BE}. Formally, in addition to the graph itself, a \LOCAL algorithm also takes as input a labeling of $G$, used to distinguish the vertices from each other.

    \begin{defn}[{\LOCAL algorithms}]\label{defn:LOCAL}
        A \emphd{\LOCAL algorithm} is a function $\mathcal{A}$ that sends isomorphism classes of finite rooted labeled graphs to $\N$. Given a locally finite labeled graph $(G,f)$ and $T \in \N$, the \emphd{output of $\mathcal{A}$ on $G$ after time $T$} is the function $\mathcal{A}_T(G, f) \colon V(G) \to \N$ given by
		\[
			\big(\mathcal{A}_T(G, f)\big) (v) \,\defeq\, \mathcal{A}\left([G, f, v]_T\right) \quad \text{for all } v \in V(G).
		\]
    \end{defn}

    Observe that an LCL problem, as defined in Definition~\ref{defn:LCL}, is a pair $\Pi = (R, \mathcal{P})$ such that $R \in \N$ and $\mathcal{P}$ is a \LOCAL algorithm taking values in $\set{0,1}$, and 
    a labeling $f \colon V(G) \to \N$ is a $\Pi$-coloring of $G$ exactly when $\mathcal{P}_R(G,f)$ is the constant $1$ function on $V(G)$.

    There are two variations of the \LOCAL model, depending on the way the symmetry-breaking labeling $f$ in Definition~\ref{defn:LOCAL} is chosen. In the deterministic version of the model, $f$ assigns a unique identifier to each vertex, and the output of the algorithm must provide a valid solution to the problem regardless of the specific assignment of the identifiers:

    \begin{defn}[Deterministic \LOCAL complexity of LCL problems]\label{defn:det_comp}
        Let $\Pi$ be an LCL problem and let $\mathbb{G}$ be a class of finite graphs. The \emphd{deterministic \LOCAL complexity} of $\Pi$ on the class $\mathbb{G}$ is the function $\mathsf{Det}_{\Pi, \mathbb{G}} \colon \N^+ \to \N \cup \set{\infty}$ defined by letting $\mathsf{Det}_{\Pi, \mathbb{G}}(n)$ be the smallest $T \in \N$ such that there exists a \LOCAL algorithm $\mathcal{A}$ with the following property:
		
		
		\begin{leftbar}
		\noindent For every $n$-vertex graph $G \in \mathbb{G}$ and every bijection $\mathsf{ID} \colon V(G) \to n$, the function $\mathcal{A}_T(G, \mathsf{ID})$ is a $\Pi$-coloring of $G$.
		\end{leftbar}

        \noindent If no such $T \in \N$ exists, we set $\mathsf{Det}_{\Pi, \mathbb{G}}(n) \defeq \infty$.
    \end{defn}

    In the randomized version of the \LOCAL model, the symmetry is broken using a random labeling:

    \begin{defn}[Randomized \LOCAL complexity of LCL problems]\label{defn:rand_comp}
        Let $\Pi$ be an LCL problem and let $\mathbb{G}$ be a class of finite graphs. The \emphd{randomized \LOCAL complexity} of $\Pi$ on the class $\mathbb{G}$ is the function $\mathsf{Rand}_{\Pi, \mathbb{G}} \colon \N^+ \to \N \cup \set{\infty}$ defined by letting $\mathsf{Det}_{\Pi, \mathbb{G}}(n)$ be the smallest $T \in \N$ such that there exist $\ell \in \N^+$ and a \LOCAL algorithm $\mathcal{A}$ with the following property:
		
		
		\begin{leftbar}
		\noindent For every $n$-vertex graph $G \in \mathbb{G}$, if a mapping $\theta \colon V(G) \to \ell$ is chosen uniformly at random, then
        \[
            \P\big[\text{the function $\mathcal{A}_T(G, \theta)$ is a $\Pi$-coloring of $G$}\big] \,\geq\, 1 - \frac{1}{n}.
        \]
		\end{leftbar}

        \noindent If no such $T \in \N$ exists, we set $\mathsf{Det}_{\Pi, \mathbb{G}}(n) \defeq \infty$.
    \end{defn}

    Unsurprisingly, the randomized version of the model is more computationally powerful than the deterministic one, and there exist many instances where a problem's randomized \LOCAL complexity is significantly lower than its deterministic \LOCAL complexity \cite{CKP}.

    In \cite{BernshteynDistributed}, the first named author showed that fast deterministic \LOCAL algorithms for LCL problems yield Borel solutions to these problems:\footnote{For simplicity, Theorems \ref{theo:dist_Borel} and \ref{theo:dist_meas} are presented here in a somewhat less general way than the corresponding statements in \cite{BernshteynDistributed}. However, these less general formulations are sufficient for most applications. See the discussion in \S\ref{subsec:LOCAL_proof} and \cite[\S2.B]{BernshteynDistributed} for more details.}

    \begin{theo}[{AB \cite[Theorem 2.10]{BernshteynDistributed}}]\label{theo:dist_Borel}
        Let $\Pi$ be an LCL problem and let $\mathbb{G}$ be a class of $\Pi$-colorable finite graphs such that $\mathsf{Det}_{\Pi,\mathbb{G}}(n) = o(\log n)$. If $G$ is a Borel graph of finite maximum degree all of whose finite induced subgraphs are in $\mathbb{G}$, then $G$ has a Borel $\Pi$-coloring.
    \end{theo}

    Similarly, a fast {randomized} \LOCAL algorithm yields measurable/Baire-measurable solutions:

    \begin{theo}[{AB \cite[Theorem 2.14]{BernshteynDistributed}}]\label{theo:dist_meas}
        Let $\Pi$ be an LCL problem and let $\mathbb{G}$ be a class of $\Pi$-colorable finite graphs such that $\mathsf{Rand}_{\Pi,\mathbb{G}}(n) = o(\log n)$. If $G$ is a Borel graph of finite maximum degree all of whose finite induced subgraphs are in $\mathbb{G}$, then the following conclusions hold:
        \begin{enumerate}[label=\ep{\normalfont\roman*}]
			\item\label{item:LLLmeas} If $\mu$ is a probability Borel measure on $V(G)$, then $G$ has a $\mu$-measurable $\Pi$-coloring.
			
			\item\label{item:LLLBaire} If $\tau$ is a compatible Polish topology on $V(G)$, then $G$ has a $\tau$-Baire-measurable $\Pi$-coloring. 
		\end{enumerate}
    \end{theo}

    In other words, if we are willing to sacrifice Borelness and settle for measurable or Baire-measurable solutions, then we may use the much more powerful randomized version of the \LOCAL model to find them. 
    Here we show that if the graph $G$ has finite asymptotic separation index, then no sacrifice is necessary: we can obtain \emph{Borel} $\Pi$-colorings of $G$ from \emph{randomized} \LOCAL algorithms:

    \begin{tcolorbox}  
    \begin{theo}[Borel colorings using randomized algorithms and finite $\asi$]\label{theo:rand_LOCAL_to_asi}
        Let $\Pi$ be an LCL problem and let $\mathbb{G}$ be a class of $\Pi$-colorable finite graphs such that $\mathsf{Rand}_{\Pi,\mathbb{G}}(n) = o(\log n)$. Let $G$ be a Borel graph of finite maximum degree all of whose finite induced subgraphs are in $\mathbb{G}$. If $\asi(G) < \infty$, then $G$ has a Borel $\Pi$-coloring.
    \end{theo}
    \end{tcolorbox}

    We in fact prove somewhat stronger (and numerically explicit) versions of Theorem~\ref{theo:rand_LOCAL_to_asi}, which are presented in \S\ref{subsec:LOCAL_proof}. We should also note that Theorems~\ref{theo:dist_Borel} and \ref{theo:dist_meas} are stated in \cite{BernshteynDistributed} in the more general context of so-called structured graphs, i.e., graphs equipped with additional combinatorial data such as an orientation, a labeling, an ordering of the vertices, etc. Our results also hold for structured graphs, as we explain in \S\ref{subsec:LOCAL_proof}.

    Although it may at first seem unrelated, Theorem~\ref{theo:rand_LOCAL_to_asi} is actually a relatively straightforward consequence of our main results concerning Borel versions of the LLL. As demonstrated in \cite{BernshteynDistributed}, randomized \LOCAL algorithms and the LLL are closely linked. Specifically, \cite[Lemma 4.6]{BernshteynDistributed} gives a way to reduce an LCL problem $\Pi$ to an application of the LLL with a \pw criterion using a sufficiently fast randomized \LOCAL algorithm. As a result, Theorem~\ref{theo:rand_LOCAL_to_asi} follows in a fairly routine manner from Theorem~\ref{theo:main_easy_asi} and \cite[Lemma 4.6]{BernshteynDistributed}; the details are presented in \S\ref{subsec:LOCAL_proof}. Curiously, the relationship between \LOCAL algorithms and the LLL is reciprocal: in order to prove Theorem~\ref{theo:main_speedup}, we combine Theorem~\ref{theo:rand_LOCAL_to_asi} with a fast randomized \LOCAL algorithm for solving instances of the LLL due to Fischer and Ghaffari \cite{FG} and Ghaffari, Harris, and Kuhn \cite{GHK}; see \S\ref{subsec:speedup_proof} for the details.

    \subsection{Combinatorial applications}\label{subsec:intro_appl}

    \subsubsection{Bounds on the Borel chromatic number}
    
    We have already mentioned two applications of our general results to the study of Borel proper colorings of graphs with finite asymptotic separation index: Theorem~\ref{theo:BorelBrooks} (Borel version of Brooks's theorem) and Theorem~\ref{theo:BorelJohansson} (Borel version of Johansson's theorem). Both of these results are derived from Theorem~\ref{theo:rand_LOCAL_to_asi} by citing relevant randomized \LOCAL algorithms.

    \begin{scproof}[ of Theorem~\ref{theo:BorelBrooks}]
        Fix an integer $\Delta \geq 3$ and let $\mathbb{B}\mathrm{rooks}_\Delta$ be the class of all finite graphs of maximum degree at most $\Delta$ with no $(\Delta+1)$-cliques. Ghaffari, Hirvonen, Kuhn, and Maus \cite{GHKM} developed a randomized \LOCAL algorithm for $\Delta$-coloring $n$-vertex graphs in the class $\mathbb{B}\mathrm{rooks}_\Delta$ with running time $O((\log \log n)^2) = o(\log n)$. If $G$ is a Borel graph of maximum degree $\Delta$ and without a $(\Delta+1)$-clique, then all its finite induced subgraphs are in $\mathbb{B}\mathrm{rooks}_\Delta$, so if $\asi(G) < \infty$, then $G$ has a Borel proper $\Delta$-coloring by Theorem~\ref{theo:rand_LOCAL_to_asi}, as desired.
    \end{scproof}

    Theorem~\ref{theo:BorelJohansson} is proved in exactly the same way by invoking randomized \LOCAL algorithms due to Chung, Pettie, and Su \cite{CPS}; see \S\ref{subsec:Joh_proof} for the details.
    As one more illustration, we show that if the Borel chromatic number of a Borel graph $G$ with finite asymptotic separation index is close to the maximum degree of $G$, then it must be equal to the ordinary chromatic number of $G$:

    \begin{tcolorbox}
    \begin{theo}[Borel colorings with $\Delta - O(\sqrt{\Delta})$ colors]\label{theo:BorelBrooksStrong}
        There exists a constant $\Delta_0 \in \N$ with the following property. Let $G$ be a Borel graph with finite maximum degree $\Delta \geq \Delta_0$ and let $k \defeq \lfloor \sqrt{\Delta + 1/4} - 7/2 \rfloor$. If $\asi(G) < \infty$ and $\chi_\mathsf{B} (G) \geq \Delta - k$, then $\chi_\mathsf{B}(G) = \chi(G)$.
    \end{theo}
    \end{tcolorbox}
    \begin{scproof}
        Clearly, $\chi(G) \leq \chi_\mathsf{B}(G)$. To show $\chi(G) \geq \chi_\mathsf{B}(G)$, set $q \defeq \chi_\mathsf{B}(G) - 1$ and suppose, toward a contradiction, that $\chi(G) \leq q$. 
        Let $\Pi_q$ be the LCL problem that encodes proper $q$-coloring of graphs \ep{see Example~\ref{exmp:coloringLCL}} and let $\mathbb{C}\mathrm{ol}_{\Delta, q}$ be the class of all finite graphs with maximum degree at most $\Delta$ and chromatic number at most $q$. For large $\Delta$ and $q \geq \Delta - k - 1$, Bamas and Esperet \cite[Theorem 1.3]{BamasEsperet} proved that
		\[
			\mathsf{Rand}_{\Pi_{q}, \mathbb{C}\mathrm{ol}_{\Delta, q}}(n) \,\leq\, \exp(O(\sqrt{\log \log n})) \,=\, o(\log n).
		\]
        Since all finite induced subgraphs of $G$ are in $\mathbb{C}\mathrm{ol}_{\Delta, q}$, Theorem~\ref{theo:rand_LOCAL_to_asi} implies that $G$ has a Borel proper $q$-coloring, which contradicts the choice of $q$. 
    \end{scproof}

    When $\Delta$ is sufficiently large, Theorem~\ref{theo:BorelBrooks} is a consequence of Theorem~\ref{theo:BorelBrooksStrong} and Brooks's Theorem~\ref{theo:Brooks} for the ordinary chromatic number. Similarly, we have:

    \begin{tcolorbox}
    \begin{corl}[Borel colorings with $\Delta - 1$ colors]\label{theo:BorelReed}
        There exists a constant $\Delta_0 \in \N$ with the following property. Let $G$ be a Borel graph with finite maximum degree $\Delta \geq \Delta_0$ and without a $\Delta$-clique. If $\asi(G) < \infty$, then $\chi_\mathsf{B}(G) \leq \Delta - 1$.
    \end{corl}
    \end{tcolorbox}
    \begin{scproof}
        For large $\Delta$, Reed \cite{ReedDelta1} showed that graphs with maximum degree $\Delta$ and without $\Delta$-cliques have chromatic number at most $\Delta - 1$. Hence, the desired result follows from Theorem~\ref{theo:BorelBrooksStrong}.
    \end{scproof}

    \subsubsection{Borel edge-colorings}\label{subsubsec:edge-coloring}

    Alongside proper coloring, another classical concept in graph theory is proper edge-coloring, which is defined as follows:

    \begin{defn}[Edge-colorings and the \ep{Borel} chromatic index]
        Let $G$ be a graph and let $q \in \N$. A \emphd{proper $q$-edge-coloring} of $G$ is a function $f \colon E(G) \to q$ such that $f(e) \neq f(e')$ for all distinct edges $e$, $e' \in E(G)$ that share a vertex. The \emphd{chromatic index} $\chi'(G)$ of $G$ is the smallest $q \in \N$ such that $G$ has a proper $q$-edge-coloring; if no such $q \in \N$ exists, we set $\chi'(G) \defeq \infty$. Similarly, if $G$ is a Borel graph, then the \emphd{Borel chromatic index} $\chi'_\mathsf{B}(G)$ of $G$ is the smallest $q \in \N$ such that $G$ has a Borel proper $q$-edge-coloring; 
        if no such $q$ exists, we set $\chi'_\mathsf{B}(G) \defeq \infty$.
    \end{defn}

    A proper $q$-edge-coloring of a graph $G$ is exactly the same thing as a proper $q$-coloring of the \emphd{line graph} of $G$, i.e., the graph $L(G)$ with vertex set $E(G)$ and edge set $\set{\set{e,e'} \in [E(G)]^2 \,:\, |e \cap e'| = 1}$. If $G$ is a Borel graph of finite maximum degree $\Delta$, then the maximum degree of $L(G)$ is at most $2\Delta -2$, and hence the Borel chromatic index of $G$ is at most $2\Delta - 1$ by the Kechris--Solecki--Todorcevic Theorem~\ref{theo:KST}. Marks \cite{Marks} proved that there exist Borel forests $G$ for which this bound is optimal, and Conley, Jackson, Marks, Seward, and Tucker-Drob \cite{HypComb} further showed that such $G$ can be taken to be hyperfinite. The following theorem summarizes these facts:

    \begin{theo}[{Kechris--Solecki--Todorcevic \cite[15]{KST}, Marks \cite[Theorem 1.4]{Marks}, Conley--Jackson--Marks--Seward--Tucker-Drob \cite[Theorem 1.4(2)]{HypComb}}]\label{theo:Borel_edge_coloring}
        Fix an integer $\Delta \in \N$.

        \begin{enumerate}[label=\ep{\normalfont\roman*}]
            \item\label{item:upper} Every Borel graph $G$ of maximum degree $\Delta$ satisfies $\chi'_\mathsf{B}(G) \leq 2\Delta - 1$.

            \item\label{item:lower} There exists a hyperfinite Borel forest $G$ of maximum degree $\Delta$ with $\chi'_\mathsf{B}(G) = 2\Delta - 1$.
        \end{enumerate}
    \end{theo}

    By contrast, a celebrated result of Vizing asserts that the ordinary chromatic index of a graph of maximum degree $\Delta$ cannot exceed $\Delta + 1$: 

    \begin{theo}[{Vizing \cite{Vizing}}]\label{theo:Vizing}
        If $G$ is a graph of finite maximum degree $\Delta$, then $\chi'(G) \leq \Delta + 1$. 
    \end{theo}
    
    See \cite[\S{}A.1]{EdgeColoringMonograph} for an English translation of Vizing's paper and \cites[\S17.2]{BondyMurty}[\S5.3]{Diestel} for modern textbook presentations.
    
    Theorem~\ref{theo:Borel_edge_coloring} precludes a Borel version of Vizing's theorem when $\Delta \geq 3$. On the other hand, a measurable version of Vizing's theorem was recently proved by Greb\'ik \cite{GmeasVizing} (building on the earlier breakthrough work of Greb\'ik and Pikhurko \cite{GP}). 
%
   %
    Naturally, we would like to know what additional assumptions on $G$ bring the Borel chromatic index of $G$ closer to $\chi'(G)$. For example, Dhawan and the first named author \cite{SubexpVizing} showed that $\chi'_\mathsf{B}(G) \leq \Delta + 1$ if $G$ is of subexponential growth. In this direction, the asymptotic separation index turns out to be an extremely helpful parameter, as demonstrated by Qian and the second named author \cite{ASIalgorithms}:

    \begin{theo}[{Qian--FW \cite[Theorem 3.7]{ASIalgorithms}}]\label{theo:Borel_edge-coloring_asi}
        If $G$ is a Borel graph with finite maximum degree $\Delta$ and with $\asi(G) < \infty$, then $\chi_\mathsf{B}'(G) \leq \chi'(G) + \asi(G)$, and hence $\chi_\mathsf{B}'(G) \leq \Delta + 1 + \asi(G)$.
    \end{theo}

    In the case when $G$ has no odd cycles, Theorem~\ref{theo:Borel_edge-coloring_asi} was established earlier by Bowen and the second named author \cite{BWKonig}. Note that for graphs $G$ with $\asi(G) \leq 1$, examples of which are given in \S\ref{sec:intro_asi}, Theorem~\ref{theo:Borel_edge-coloring_asi} yields the bound $\chi_\mathsf{B}'(G) \leq \Delta + 2$. It is still an open problem whether this can be reduced to $\Delta + 1$. It is also not known whether a bound of the form $\chi_\mathsf{B}'(G) \leq \Delta + c$ for some fixed $c$ (perhaps even $c = 1$) holds for all graphs $G$ with $\asi(G) < \infty$ regardless of the value of $\asi(G)$. As a consequence of our main results, we obtain an upper bound on $\chi'_\mathsf{B}(G)$ for graphs $G$ with finite asymptotic separation index that is asymptotic to $\Delta$ and independent of the value of $\asi(G)$:
    
    \begin{tcolorbox}
    \begin{theo}[Borel edge-colorings with $\approx \Delta$ colors]\label{theo:Borel_edge-coloring_close_to_Delta}
        There is $\Delta_0 \in \N$ such that if $G$ is a Borel graph of finite maximum degree $\Delta \geq \Delta_0$ and $\asi(G) < \infty$, then $\chi'_\mathsf{B}(G) \leq \Delta + \sqrt{\Delta} \log^3 \Delta.$
    \end{theo}
    \end{tcolorbox}

    Theorem~\ref{theo:Borel_edge-coloring_close_to_Delta} is proved by appealing to a randomized \LOCAL algorithm for edge-coloring due to Chang, He, Li, Pettie, and Uitto \cite{CHLPU}
    ; see \S\ref{subsec:proof_edge_epsilon} for the details. In view of Theorem~\ref{theo:Borel_edge-coloring_asi}, Theorem~\ref{theo:Borel_edge-coloring_close_to_Delta} is only interesting for graphs $G$ whose asymptotic separation index is non-negligible compared to their maximum degree. This again highlights the importance of determining whether there exist Borel graphs $G$ with $1 < \asi(G) < \infty$.

    \subsubsection{Borel edge-colorings of Schreier graphs}\label{subsubsec:Schreier}

    Some of the most well-studied examples in descriptive combinatorics are the so-called Schreier graphs \cite[\S6.5]{KechrisMarks}. These graphs are associated to group actions and arise naturally in measurable and topological dynamics. For the remainder of \S\ref{subsubsec:Schreier}, we fix a finitely generated group $\G$ with a finite symmetric generating set $F \subseteq \G \setminus \set{\mathbf{1}_\G}$, where $\mathbf{1}_\G$ is the identity element of $\G$ (a subset of a group is \emphd{symmetric} if it is closed under taking inverses). Recall that a group action $\bm{a} \colon \G \acts X$ is \emphd{free} if $\gamma \cdot_{\bm{a}} x \neq x$ for all $x \in X$ and $\gamma \in \G \setminus \set{\mathbf{1}_\G}$, i.e., if the stabilizer of every point in $X$ is trivial.

    \begin{defn}[Schreier graphs]\label{defn:Schreier}
        The \emphd{Schreier graph} of a free action $\bm{a} \colon \G \acts X$ is the graph $\mathsf{Sch}(\bm{a}, F)$  with vertex set $X$ and edge set $\set{\set{x, \sigma \cdot_{\bm{a}} x} \,:\, x \in X, \ \sigma \in F}$.
    \end{defn}

    Note that if $X$ is a standard Borel space and $\bm{a} \colon \G \acts X$ is a free Borel action (meaning that the map $x \mapsto \gamma \cdot_{\bm{a}} x$ is Borel for each $\gamma \in \G$), then the corresponding Schreier graph $\mathsf{Sch}(\bm{a},F)$ is a Borel graph of maximum degree $|F|$. In contrast to general Borel graphs, the bound in Theorem~\ref{theo:Borel_edge_coloring}\ref{item:upper} is not sharp when $G$ is a Schreier graph, as it can be easily seen \cite[67]{KechrisMarks} that 
    \[
        \chi_\mathsf{B}'(\mathsf{Sch}(\bm{a},F)) \,\leq\, |\set{\sigma \in F \,:\, \mathrm{ord}(\sigma) \text{ is even}}| \,+\, \frac{3}{2} |\set{\sigma \in F \,:\, \mathrm{ord}(\sigma) \text{ is odd or infinite}}| \,\leq\, \frac{3}{2} |F|.
    \]
    It is an open problem to determine the optimal upper bound on $\chi_\mathsf{B}'(\mathsf{Sch}(\bm{a},F))$ in terms of $|F|$.
    
    The investigation of Borel edge-colorings of Schreier graphs with finite asymptotic separation index was initiated by the second named author in \cite{FelixFinDim}. There it is shown that a Borel version of Vizing's theorem holds for Schreier graphs with asymptotic separation index at most $1$ under the assumption that no element of $F$ is of odd order (this is also a consequence of Theorem~\ref{theo:Borel_edge-coloring_asi}, which was proved later). Here we generalize this result by extending it to the case of arbitrary finite asymptotic separation index and allowing $F$ to contain elements of sufficiently large odd order:  
%
%

    \begin{tcolorbox}
    \begin{theo}[Borel edge-colorings of Schreier graphs]\label{theo:BorelVizing}
        Suppose that $|F| \geq 2$ and for each $\sigma \in F$, the order of $\sigma$ is either even or at least $2^{30}\,|F|\log|F|$. Let $\bm{a} \colon \G \acts X$ be a free Borel action of $\G$ on a standard Borel space $X$ and let $G \defeq \mathsf{Sch}(\bm{a}, F)$ be the corresponding Schreier graph. If $\asi(G) < \infty$, then $\chi'_\mathsf{B}(G) \leq |F| + 1$.
    \end{theo}
    \end{tcolorbox}


    When $|F| \leq 1$, the corresponding Schreier graphs $G$ trivially satisfy $\chi'_\mathsf{B}(G) \leq |F|$, so there is no harm in only considering the $|F| \geq 2$ case. The assumption that $|F| \geq 2$ in Theorem~\ref{theo:BorelVizing} ensures that $\log |F| > 0$. We made no attempt to optimize the constant factor $2^{30}$. Note that Theorem~\ref{theo:BorelVizing} applies to some groups of exponential growth thanks to Theorem~\ref{theo:asi1_group} (a Borel version of Vizing's theorem for graphs of subexponential growth is already known \cite{SubexpVizing}).

    Theorem~\ref{theo:BorelVizing} cannot be easily derived using Theorem~\ref{theo:rand_LOCAL_to_asi}, because no $o(\log n)$-time randomized \LOCAL algorithm for Vizing's theorem is known (although the existence of such an algorithm has not been ruled out), with the most efficient currently available algorithm running in time $O(\log^2 n)$ \cite{FastVizing}. Instead, we prove Theorem~\ref{theo:BorelVizing} by invoking our Borel LLL (i.e., Theorem~\ref{theo:main_speedup}) directly. More precisely, we apply Theorem~\ref{theo:main_speedup} to prove a certain general result about Borel independent complete sections, which is described in \S\ref{subsec:cs}, and then use it to deduce Theorem~\ref{theo:BorelVizing}.

    \subsubsection{Borel independent complete sections}\label{subsec:cs}

    For a graph $G$, a set $I \subseteq V(G)$ is \emphd{$G$-independent} if no two vertices in $I$ are adjacent. If $\mathcal{F}$ is a family of pairwise disjoint sets, then a \emphd{complete section} for $\mathcal{F}$ is a set $S$ that includes at least one element from every set in $\mathcal{F}$. The following classical fact is proved by a simple application of the LLL:

    \begin{theo}[{Alon \cites{LinArbor,AlonRef}[Proposition 5.5.3]{AS}}]\label{theo:Alon}
        Let $G$ be a graph of finite maximum degree $\Delta$ and let $\mathcal{F}$ be a family of pairwise disjoint subsets of $V(G)$. If $|S| \geq 2e\Delta$ for all $S \in \mathcal{F}$, then there is a $G$-independent complete section for $\mathcal{F}$.
    \end{theo}

    Theorem~\ref{theo:Alon}, and independent complete sections in general, have been a subject of considerable interest, not least because of their many applications to other combinatorial problems.\footnote{A related notion considered in the literature is that of an independent transversal---i.e., a set that meets each member of $\mathcal{F}$ in \emph{exactly} one point. For our purposes, working with complete sections will be more natural. Note that the existence of a $G$-independent complete section for $\mathcal{F}$ is equivalent to the existence of a $G$-independent transversal (although this equivalence may fail in the Borel setting).} The required lower bound on the sizes of the sets in $\mathcal{F}$ has been subsequently improved from $2e\Delta$ to $2\Delta$ by Haxell \cite{Haxell}, which is best possible \cite{BES}.

    We prove a Borel version of Theorem~\ref{theo:Alon} under a finite asymptotic separation index assumption. To simplify the statement of our result, we say that a set $S \subseteq V(G)$ is a \emphd{complete section} for a graph $G$ if it meets every connected component of $G$. In other words, $S$ is a complete section for $G$ if and only if it is a complete section for the following family of sets:
    \[
        \mathcal{F} \,\defeq\, \set{V(C) \,:\, \text{$C$ is a connected component of $G$}}.
    \]
    A subgraph $G'$ of $G$ is \emphd{spanning} if $V(G') = V(G)$.
    

    \begin{tcolorbox}
    \begin{theo}[Borel independent complete sections]\label{theo:indep_transversal}
        Let $H$ be a Borel graph with $\asi(H) < \infty$ and let $G_1$, $G_2$ be two spanning Borel subgraphs of $H$. Suppose that the maximum degree of $G_1$ is at most $\Delta \in \N$, where $\Delta \geq 2$. 
        If every component of $G_2$ contains at least $
            2^{25} \,\Delta\log\Delta
        $ 
        vertices, then there is a Borel $G_1$-independent complete section for $G_2$.
    \end{theo}
    \end{tcolorbox}

    The assumption that $\asi(H) < \infty$ in Theorem~\ref{theo:indep_transversal} cannot be removed. Indeed, Marks \cite[Theorem 1.6]{Marks} gave an example of a Borel graph $H$ and two spanning Borel subgraphs $G_1$, $G_2$ such that both $G_1$ and $G_2$ have maximum degree $2$, all their components are infinite, and every Borel set $A \subseteq V(H)$ either includes an entire component of $G_1$ (and hence it is definitely not $G_1$-independent), or misses a component of $G_2$ (i.e., it is not a complete section for $G_2$).

    The usual proof of Theorem~\ref{theo:Alon} is one of the rare applications of the LLL for which the bound $\pr(\de+1) \leq 1/e$ is sharp and a \pw criterion is insufficient. As a result, we cannot establish Theorem~\ref{theo:indep_transversal} by simply following the standard proof of Theorem~\ref{theo:Alon} and substituting Theorem~\ref{theo:main_easy_asi} or Theorem~\ref{theo:main_speedup} in place of the LLL. Nevertheless, we are able to prove Theorem~\ref{theo:indep_transversal} using a different approach, which unfortunately results in a superlinear lower bound on the size of the components of $G_2$. However, we suspect that the conclusion of Theorem~\ref{theo:indep_transversal} should hold even if the components of $G_2$ have size bounded below by some linear function of $\Delta$:

    \begin{conj}\label{conj:cs}
        There exist constants $C$, $\Delta_0 > 0$ with the following property: Let $H$ be a Borel graph with $\asi(H) < \infty$ and let $G_1$, $G_2$ be two spanning Borel subgraphs of $H$. Suppose that the maximum degree of $G_1$ is at most $\Delta \in \N$, where $\Delta \geq \Delta_0$. 
        If every component of $G_2$ contains at least $
            C\Delta
        $ 
        vertices, then there is a Borel $G_1$-independent complete section for $G_2$.
    \end{conj}

     One way to settle Conjecture~\ref{conj:cs} would be to prove a version of Theorem~\ref{theo:main_speedup} under the usual LLL assumption $\pr(\de+1) \leq 1/e$. Since we suspect that $\Delta \log \Delta$ is not the optimal order of magnitude in Theorem~\ref{theo:indep_transversal}, we made no attempt to optimize the constant factor $2^{25}$. We prove Theorem~\ref{theo:indep_transversal} in \S\ref{subsec:proof_cs} and then use it to derive Theorem~\ref{theo:BorelVizing} in \S\ref{subsec:proof_Schreier}. If Conjecture~\ref{conj:cs} is true, then the same argument would improve Theorem~\ref{theo:BorelVizing} by replacing the bound $2^{30}\,|F|\log|F|$ by a linear function of $|F|$. We expect that this can also be achieved by following a ``nibbling'' approach similar to the one in the paper \cite{ReedSud} by Reed and Sudakov, but the resulting argument would be quite technical, and since we do not see how to generalize it to prove Conjecture~\ref{conj:cs} in full, we do not pursue it in this paper.


    \section{Preliminaries}

    \subsection{The spaces of finite sets and constraints}\label{subsec:const}

    In this subsection, given a standard Borel space $X$ and an integer $q \in \N^+$, we describe a way of endowing the set $\Const(X,q)$ of all $(X,q)$-constraints with the structure of a standard Borel space (this is necessary for Definition~\ref{defn:BorelCSP} to make sense). The material in this subsection is entirely standard; nevertheless, we decided to include it here for completeness.
    
    Let $X$ be a standard Borel space. We begin by equipping $[X]^k$ for $k \in \N$ with the structure of a standard Borel space. To this end, let $(X)^k$ be the set of all $k$-tuples $(x_0, \ldots, x_{k-1}) \in X^k$ of distinct points, and let $\sim_k$ be the equivalence relation on $(X)^k$ given by
    \[
        (x_0, \ldots, x_{k-1}) \,\sim_k\, (y_0, \ldots, y_{k-1}) \quad \vcentcolon\Longleftrightarrow \quad \set{x_0, \ldots, x_{k-1}} \,=\, \set{y_0, \ldots, y_{k-1}}.
    \]
    By \cite[Example 6.1 and Proposition 6.3]{KechrisMiller}, the quotient space $[X]^k = (X)^k/\sim_k$ is standard Borel. It follows that the disjoint union $\finset{X} = [X]^0 \sqcup [X]^1 \sqcup [X]^2 \sqcup \ldots$ is standard Borel as well \cite[Proposition 1.4]{AnushDST}.
    
    Next, given $q \in \N^+$, we let $\finfun{X}{q}$ be the set of all partial functions $\phi \colon X \pto q$ such that $\dom(\phi)$ is a finite subset of $X$. Each $\phi \in \finfun{X}{q}$ can be identified with a finite subset of $X \times q$, namely, with its graph $\set{(x, \phi(x)) \,:\, x \in \dom(\phi)}$. Using this identification, it is routine to check that $\finfun{X}{q}$ is a Borel subset of $\finset{X \times q}$.
    
    Finally, since an $(X,q)$-constraint is a finite subset of $\finfun{X}{q}$, we may define
    \[
        \Const(X,q) \,\defeq\, \big\{B \in \finset{\finfun{X}{q}} \,:\, \text{all functions in $B$ have the same domain} \big\}.
    \]
    It can be easily verified that $\Const(X,q)$ is a Borel subset of the standard Borel space $\finset{\finfun{X}{q}}$, and therefore, equipped with the relative $\sigma$-algebra of Borel sets, it is also standard Borel.

    \subsection{Borel combinatorics on component-finite graphs}\label{subsec:cf}

    As mentioned in the introduction, the reason we expect graphs with finite asymptotic separation index to have well-behaved Borel combinatorics is that Borel combinatorics ``trivialize'' on component-finite graphs. The following statement makes this intuition rigorous:

    \begin{prop}[Borel combinatorics on component-finite graphs]\label{prop:cf}
        Let $\B \colon X \to^? q$ be a Borel CSP on a standard Borel space $X$. If the graph $G_\B$ is component-finite and $\B$ has a solution, then it has a Borel solution as well.
    \end{prop}

    The proof of Proposition~\ref{prop:cf} relies on the following deep and classical result (which we shall use throughout the paper without mention):
    
    \begin{theo}[{Luzin--Novikov uniformization \cite[Theorem 18.10]{KechrisDST}}]\label{theo:LN}
        Let $X$, $Y$ be standard Borel spaces and let $R \subseteq X \times Y$ be a Borel subset. If for every $x \in X$, the set
        \[
            R_x \,\defeq\, \set{y \in Y \,:\, (x,y) \in R}
        \]
        is countable, then there exists a sequence $(f_n)_{n \in \N}$ of Borel partial functions $f_n \colon X \pto Y$ defined on Borel subsets of $X$ such that for all $x \in X$, $R_x = \set{f_n(x) \,:\, n \in \N,\ x \in \dom(f_n)}$.
        
        If additionally each set $R_x$ is nonempty, then the mapping $f \colon X \to Y$ given by
        \[
            f(x) \,\defeq\, f_{n(x)}(x), \text{ where } n(x) \,\defeq\, \min \set{n \in \N \,:\, x \in \dom(f_n)},
        \]
        is a Borel function such that $f(x) \in R_x$ for all $x \in X$.
    \end{theo}



    \begin{scproof}[ of Proposition~\ref{prop:cf}]
        This is a standard argument in descriptive combinatorics, so we only provide a proof sketch. Roughly, since on each component of $G_\B$ the problem can be solved in only finitely many ways, the Luzin--Novikov theorem yields a Borel choice of one solution on every component; putting the chosen solutions together solves the problem on the whole space.
        
        In a little more detail, let $\mathcal{C} \defeq \set{V(C) \,:\, \text{$C$ is a connected component of $G_\B$}}$. This is a Borel subset of $\finset{X}$. Say that a \emphd{finite partial solution} to $\B$ is a map $\phi \in \finfun{X}{q}$ such that $\phi$ satisfies every constraint $B \in \B$ with $\dom(B) \subseteq \dom(\phi)$. Then the set $\mathcal{S} \subseteq \finfun{X}{q}$ of all finite partial solutions to $\B$ is also Borel. Since $\B$ has a solution, for each finite set $U \subseteq X$, there is some $\phi \in \mathcal{S}$ with $\dom(\phi) = U$. Let 
        \[
            R \,\defeq\, \set{(U, \phi) \in \mathcal{C} \times \mathcal{S} \,:\, \dom(\phi) = U}.
        \]
        Then $R$ is a Borel subset of $\mathcal{C} \times \mathcal{S}$, and for each $U \in \mathcal{C}$, the set $R_U \defeq \set{\phi \in \mathcal{S} \,:\, (U, \phi) \in R}$ is finite and nonempty. Therefore, by the Luzin--Novikov theorem (Theorem~\ref{theo:LN}), there exists a Borel function $\mathcal{C} \to \mathcal{S} \colon U \mapsto \phi_U$ such that $\dom(\phi_U) = U$ for all $U$. For each $x \in X$, let $[x] \in \mathcal{C}$ be the vertex set of the component containing $x$. Then the map $f \colon X \to q$ given by
        \[
            f(x) \,\defeq\, \phi_{[x]}(x) \quad \text{for all } x \in X,
        \]
        or, equivalently, $f \defeq \bigsqcup_{U \in \mathcal{C}} \phi_U$, is a Borel solution to $\B$, as desired.
    \end{scproof}

    \begin{corl}\label{corl:coloring}
        If $G$ is a component-finite Borel graph and $\chi(G) < \infty$, then $\chi_\mathsf{B}(G) = \chi(G)$.
    \end{corl}
    \begin{scproof}
        Clearly, $\chi_\mathsf{B}(G) \geq \chi(G)$. To see that $\chi_\mathsf{B}(G) \leq \chi(G)$, recall from Example \ref{exmp:coloringBorelCSP} that the proper $\chi(G)$-coloring problem for $G$ can be encoded via a Borel CSP $\B \colon V(G) \to^? \chi(G)$ given by
        \[
            \B \,\defeq\, \big\{\set{(u \mapsto i,\, v \mapsto i) \,:\, 0 \leq i < \chi(G)} \,:\, uv \in E(G)\}.
        \]
        Since $G_\B = G$ is component-finite, it has a Borel proper $\chi(G)$-coloring by Proposition~\ref{prop:cf}.
    \end{scproof}

    Another useful consequence of the Luzin--Novikov theorem is the fact that if $f \colon X \to Y$ is a countable-to-one Borel function between standard Borel spaces, then its image $f(X)$ is a Borel subset of $Y$ \cite[Exercise 18.14]{KechrisDST}.

    
    \section{Proofs of the main results}
    
    \subsection{Proof of Theorem \ref{theo:main_easy_asi}}\label{subsec:proof_main_easy}

    As stated in \S\ref{subsec:intro_Borel_LLL}, in our proof of Theorem~\ref{theo:main_easy_asi} we will use an auxiliary parameter in place of the separation index. The following definition is reminiscent of (and inspired by) an alternative way of defining asymptotic separation index given in \cite[Lemma 3.1(1')]{Dimension}.

    \begin{defn}[Shattering number]\label{defn:shatter}
        Let $\B \colon X \to^? q$ be a Borel CSP on a standard Borel space $X$. The \emphd{shattering number} of $\B$, denoted by $\mathsf{sh}(\B)$, is the smallest $s \in \N$ for which there exists a Borel family $\mathcal{P} \subseteq \finset{X}$ of finite subsets of $X$ such that $\mathcal{P}$ is a partition of $X$ and for each $B \in \B$, $\dom(B)$ intersects at most $s$ members of $\mathcal{P}$ (if no such $s \in \N$ exists, we set $\mathsf{sh}(\B) \defeq \infty$).
    \end{defn}

    Some readers may be more familiar with describing a partition of a standard Borel space $X$ into finite sets as an equivalence relation on $X$ with finite classes. For their benefit, we remark that a partition $\mathcal{P}$ is a Borel subset of $\finset{X}$ if and only if the equivalence relation whose classes are the members of $\mathcal{P}$ is a Borel subset of $X^2$ \cite[Lemma 3.8]{BernshteynYu}.
    
    Let us observe two easy upper bounds on the shattering number of a Borel CSP.

    \begin{lemma}\label{lem:sh}
        If $\B \colon X \to^? q$ be a Borel CSP on a standard Borel space $X$, then:
        \begin{enumerate}[label=\ep{\normalfont\roman*}]
            \item\label{item:si} $\mathsf{sh}(\B) \leq \si(G_\B) + 1$, and
            \item\label{item:dom} $\mathsf{sh}(\B) \leq \sup_{B \in \B} |\dom(B)|$.
        \end{enumerate}
    \end{lemma}
    \begin{scproof}
        For \ref{item:si}, suppose that $X = U_0 \sqcup \ldots \sqcup U_s$ is a partition witnessing $\si(G_\B) = s$. Then
    \[
        \mathcal{P} \,\defeq\, \set{V(C) \,:\, \text{$C$ is a connected component of $G[U_i]$ for some $0 \leq i \leq s$}}
    \]
    is a family of finite sets witnessing $\mathsf{sh}(\B) \leq s + 1$. For \ref{item:dom}, take $\mathcal{P} = \set{\set{x} \,:\, x \in X}$.
    \end{scproof}


    We can now state the strengthening of Theorem~\ref{theo:main_easy_asi} using the shattering number: 

    \begin{tcolorbox}
    \begin{theo}[Borel LLL with bounded shattering number]\label{theo:main_easy_shatter}
        Let $\B \colon X \to^? q$ be a Borel CSP on a standard Borel space $X$. Suppose that $\mathsf{sh}(\B) \leq s$ for some $s \in \N^+$. If
        \begin{equation}\label{eq:shatterbound}
            \pr(\B) \, (\de(\B) + 1)^{s} \,\leq\, e^{-s},
        \end{equation}
        then $\B$ has a Borel solution $f \colon X \to q$. 
    \end{theo}
    \end{tcolorbox}

    Theorem~\ref{theo:main_easy_asi} is a direct consequence of Theorem~\ref{theo:main_easy_shatter} and Lemma~\ref{lem:sh}\ref{item:si}. On the other hand, combining Theorem~\ref{theo:main_easy_shatter} with Lemma \ref{lem:sh}\ref{item:dom} gives a conclusion reminiscent of, though slightly weaker than, \cite[Corollary 1.8]{Ber_cont}. 
    
    The technique we employ to prove Theorem~\ref{theo:main_easy_shatter} is known as the \emphd{method of conditional probabilities} in computer science, where it is commonly used for derandomization \cites[\S16]{AS}[\S5.6]{RandAlg}. This method has already been applied to obtain ``constructive'' versions of the LLL, for instance by Beck \cite{Beck}, Fischer and Ghaffari \cite{FG}, and the first named author \cite{Ber_cont}. We especially highlight the connection to \cite[Theorem 3.6]{FG} by Fischer and Ghaffari, which uses so-called network decompositions to design efficient distributed algorithms in a way that is rather analogous to our use of the shattering number.  

    Our proof of Theorem~\ref{theo:main_easy_shatter} is inductive: we shall construct the desired $q$-coloring of $X$ that solves $\B$ by iteratively coloring subsets of $X$. Throughout the iterations, the domains of the constraints will be shrinking until they become empty. Note that there are precisely two constraints with the empty domain: $\0$ and $\set{\0}$. The constraint $\set{\0}$ has probability $1/q^0 = 1$ and is violated by every coloring, while the constraint $\0$ has probability $0$ and is always satisfied.  The distinction between these two cases will play a decisive role in the last stage of our proof.
    
    To facilitate the analysis of the iterative process, the following notation will be useful. Fix a set $X$ and an integer $q \in \N^+$ and let $f \colon U \to q$ be a $q$-coloring of a subset $U \subseteq X$. For an $(X,q)$-constraint $B$, we write $B/f$ to denote the constraint with $\dom(B/f) \defeq \dom(B) \setminus U$ given by
	\[
	   B/f \,\defeq\, \set{\phi \colon (\dom(B) \setminus U) \to q \,:\, \rest{f}{\dom(B) \cap U} \,\sqcup\, \phi \,\in\, B}.
	\]
	In other words, $\phi \in B/f$ if and only if the coloring $f \sqcup \phi$ violates $B$. We emphasize that it is possible to have $\dom(B) \subseteq U$, in which case $\dom(B/f) = \0$. More specifically, if $\dom(B) \subseteq U$, then 
    \[
        B/f \,=\, \begin{cases}
            \set{\0} &\text{if $f$ violates $B$};\\
            \0 &\text{if $f$ satisfies $B$}.
        \end{cases}
    \]
    Observe that if $f' \colon U' \to q$ is another $q$-coloring of some set $U' \subseteq X \setminus U$, then \[(B/f)/f' \,=\, B/(f \sqcup f').\] Given a CSP $\B \colon X \to^? q$, we define a CSP $\B/f \colon (X \setminus U) \to^? q$ by \[\B/f \,\defeq\, \set{B/f \,:\, B \in \B}.\]
    By construction, $f' \colon (X \setminus U) \to q$ is a solution to $\B/f$ if and only if $f \sqcup f'$ is a solution to $\B$. Note that we always have $\de(\B/f) \leq \de(\B)$, but there is no \emph{a priori} relation between $\pr(\B/f)$ and $\pr(\B)$.
    
    The following lemma describes the inductive step in the proof of Theorem~\ref{theo:main_easy_shatter}:

    \begin{lemma}\label{lemma:step}
        Let $\B \colon X \to^? q$ be a Borel CSP on a standard Borel space $X$ and let $s \colon \B \to \N$ be a Borel function such that for all $B \in \B$,
        \begin{equation}\label{eq:key_lemma_assumption}
            \P[B] \, (\de(\B) + 1)^{s(B)} \,<\, e^{-s(B)}.
        \end{equation}
        Let $U \subseteq X$ be a Borel set such that the graph $G_\B[U]$ is component-finite. For each $B \in \B$, define
        \[
            \eta(B) \,\defeq\, \begin{cases}
                1 &\text{if $\dom(B) \cap U \neq \0$};\\
                0 &\text{otherwise}.
            \end{cases}
        \]
        Then there exists a Borel function $f \colon U \to q$ such that for all $B \in \B$,
        \begin{equation}\label{eq:key_lemma_conclusion}
            \P[B/f] \, (\de(\B/f) + 1)^{s(B) - \eta(B)} \,<\, e^{-s(B) + \eta(B)}.
        \end{equation}
    \end{lemma}
    \begin{scproof}
        Note that if $B \in \B$ is a constraint such that $\dom(B) \cap U = \0$, then \eqref{eq:key_lemma_conclusion} holds for this $B$ regardless of the choice of $f \colon U \to q$, since $B/f = B$, $\de(\B/f) \leq \de(\B)$, and $\eta(B) = 0$. Therefore, we only need to focus on the constraints $B \in \B$ with $\dom(B) \cap U \neq \0$. Consider any such constraint $B$. For a function $\psi \colon (\dom(B) \cap U) \to q$, let the \emphd{probability of $B$ conditioned on $\psi$} be
        \[
            \P[B \,\vert\, \psi] \,\defeq\, \frac{|\set{\phi \in B \,:\, \rest{\phi}{\dom(B) \cap U} = \psi}|}{q^{|\dom(B) \setminus U|}}.
        \]
        In other words, $\P[B \,\vert\, \psi]$ is the probability that a uniformly random $q$-coloring of $\dom(B)$ violates $B$ given that it agrees with $\psi$ on $\dom(B) \cap U$. Observe that for any $f \colon U \to q$, we have
        \begin{equation}\label{eq:Bg}
            \P[B/f] \,=\, \P[B \,\vert\, \rest{f}{\dom(B) \cap U}].
        \end{equation}
        Let $B^\ast$ be the constraint with $\dom(B^\ast) \defeq \dom(B) \cap U$ given by
        \[
	   B^\ast \,\defeq\, \big\{\psi \colon (\dom(B) \cap U) \to q \ :\  \P[B \,\vert\, \psi]\, \geq\, (e (\de(\B) + 1))^{-s(B) + 1}\big\}.
	   \]
        Simple double counting shows that
        \begin{align*}
            &(e (\de(\B) + 1))^{-s(B)} \,\stackrel{\text{\eqref{eq:key_lemma_assumption}}}{>}\, \P[B] \,=\, q^{-|\dom(B) \cap U|} \sum_{\psi \colon (\dom(B) \cap U) \to q} \P[B \,\vert\, \psi] \\
            &\qquad\qquad\qquad\geq\, q^{-|\dom(B) \cap U|} \,|B^\ast|\, (e (\de(\B) + 1))^{-s(B) + 1} \,=\, \P[B^\ast] \, (e (\de(\B) + 1))^{-s(B) + 1}.
        \end{align*}
        (This is essentially Markov's inequality from probability theory.) We conclude that
        \begin{equation}\label{eq:PBstar}
            \P[B^\ast] \, (\de(\B) + 1) \,<\, 1/e.
        \end{equation}
    Now let $\B^\ast \colon U \to^? q$ be the Borel CSP defined by
        \[
            \B^\ast \,\defeq\, \set{B^\ast \,:\, B \in \B, \ \dom(B) \cap U \neq \0}.
        \]
        Since $\de(\B^\ast) \leq \de(\B)$ by construction, it follows from \eqref{eq:PBstar} that $\B$ satisfies the assumptions of the ordinary LLL, i.e., of Lemma~\ref{lemma:LLL}. Therefore, $\B^\ast$ has a solution. Note that $G_{\B^\ast} = G_\B[U]$, which is component-finite by hypothesis. Thus, $\B^\ast$ has a Borel solution $f \colon U \to q$ by Proposition~\ref{prop:cf}. We conclude from \eqref{eq:Bg}, the definition of $\B^\ast$, and the inequality $\de(\B/f) \leq \de(\B)$ that $f$ is as desired. 
    \end{scproof}

    \begin{scproof}[ of Theorem~\ref{theo:main_easy_shatter}]
        Let $\mathcal{P} \subseteq [X]^{<\infty}$ be a Borel family of finite sets witnessing that $\mathsf{sh}(\B) \leq s$. We claim there exists a partition $\mathcal{P} = \bigsqcup_{n \in \N} \mathcal{P}_n$ of $\mathcal{P}$ into countably many Borel sets such that for each $B \in \B$ and $n \in \N$, $\dom(B)$ meets at most one member of $\mathcal{P}_n$. Indeed, let $H$ be the Borel graph with vertex set $\mathcal{P}$ in which distinct sets $S$, $R \in \mathcal{P}$ are adjacent if and only if there is a constraint $B \in \B$ such that $\dom(B)$ meets both $S$ and $R$. Each set $S \in \mathcal{P}$ is finite, and every point in $S$ belongs to at most $\de(\B) < \infty$ constraints in $\B$ ($\de(\B)$ is finite by \eqref{eq:shatterbound}). Since every constraint intersects finitely many sets in $\mathcal{P}$, it follows that the graph $H$ is locally finite. Therefore, by \cite[Proposition 4.5]{KST}, there is a Borel function $c \colon \mathcal{P} \to \N$ such that $c(S) \neq c(R)$ whenever $S$ and $R$ are adjacent in $H$ (i.e., $c$ is a Borel proper $\N$-coloring of $H$). Setting $\mathcal{P}_n \defeq c^{-1}(n)$ gives the desired partition of $\mathcal{P}$.

        For $n \in \N$, let $U_n \defeq \bigcup \mathcal{P}_n$, so $X = \bigsqcup_{n \in \N} U_n$ is a partition of $X$ into Borel sets. For each $n \in \N$, each connected component of $G_\B[U_n]$ is contained in an element of $\mathcal{P}_n$, and so $G_\B[U_n]$ is component-finite. For each $B \in \B$ and $n \in \N$, let $t_n(B)$ be the number of sets $U_i$ with $i < n$ such that $\dom(B) \cap U_i \neq \0$ and set $s_n(B) \defeq s - t_n(B)$ (thus, $t_0(B) = 0$ and $s_0(B) = s$). By assumption, $t_n(B) \leq s$ for all $n \in \N$ and thus $s_n(B) \geq 0$. 
        Hypothesis \eqref{eq:shatterbound} yields 
        \[
            \P[B] \, (\de(\B) + 1)^{s_0(B)} \,\leq\, e^{-s_0(B)} \quad \text{for all } B \in \B.
        \]
        In fact, this inequality is strict, since its 
        left-hand side is rational, while its right-hand side is not. 
        We may now set $B^{(0)} \defeq B$ for all $B \in \B$ and $\B^{(0)} \defeq \B$ and then iteratively apply Lemma~\ref{lemma:step} to produce a sequence of Borel $q$-colorings $f_n \colon U_n \to q$ such that for each $n \in \N$ and $B \in \B$,
        \begin{equation}\label{eq:iterative_bound}
            \P[B^{(n+1)}] \, (\de(\B^{(n+1)}) + 1)^{s_{n+1}(B)} \,<\, e^{-s_{n+1}(B)},
        \end{equation}
        where $B^{(n+1)} \defeq B/(f_0 \sqcup \ldots \sqcup f_{n}) = B^{(n)}/f_{n}$ and $\B^{(n+1)} \defeq \B/(f_0 \sqcup \ldots \sqcup f_{n}) = \B^{(n)}/f_{n}$. We claim that the Borel function $f \defeq \bigsqcup_{n \in \N} f_n$ is a solution to $\B$. Indeed, consider any constraint $B \in \B$. To see that $f$ satisfies $B$, note that since $\dom(B)$ if finite, there is some $n \in \N$ such that \[\rest{f}{\dom(B)} \,=\, \rest{(f_0 \sqcup \ldots \sqcup f_n)}{\dom(B)}.\] Hence, by \eqref{eq:iterative_bound}, \[\P[B/f] \,=\, \P[B/(f_0 \sqcup \ldots \sqcup f_n)] \,=\, \P[B^{(n+1)}] \,<\, (e(\de(\B^{(n+1)}) + 1))^{-s_{n+1}(B)} \,\leq\, 1.\] Since $\dom(B/f) = \0$, it follows that $B/f = \0$, i.e., $f$ satisfies $B$, as desired. 
    \end{scproof}

    \subsection{From the LLL to \LOCAL algorithms}\label{subsec:LOCAL_proof}

    
    As briefly mentioned in the introduction, LCL problems and \LOCAL algorithms can be defined in the situation when the underlying graph is equipped with additional combinatorial structure, such as an orientation, a labeling, etc. Since we will need this general set-up for the proof of Theorem~\ref{theo:main_speedup}, we recall the necessary definitions here. Our presentation follows that in \cite[\S2.A]{BernshteynDistributed}.

    For a set $X$, let $X^{< \infty} \defeq X^0 \sqcup X^1 \sqcup X^2 \sqcup \ldots$ be the set of all finite tuples of elements of $X$.
    
    \begin{defn}[{Structured graphs}]\label{defn:str_graph}
		Let $G$ be a graph. A \emphd{structure map} on $G$ is a partial function $\sigma \colon V(G)^{<\infty} \pto \N$ such that for some $\ell \in \N$, every tuple $\bm{x} \in \dom(\sigma)$ is of length at most $\ell$. A \emphd{structured graph} is a pair $\bm{G} = (G, \sigma)$, where $G$ is a graph and $\sigma$ is a structure map on $G$.
	\end{defn}

    \begin{exmps}
        \begin{enumerate}[wide, label=\ep{\normalfont\roman*}]
            \item A labeled graph $(G,f)$ in the sense of Definition~\ref{defn:LCL} can be viewed as a structured graph by identifying $V(G)$ with $V(G)^1$ in the obvious way.

            \item A directed graph $G$ can be viewed as a structured graph by encoding the directions of the arcs of $G$ using the function $\sigma \colon V(G)^2 \to 2$ that sends a pair $(u, v)$ to $1$ if and only if there is an arc from $u$ to $v$ in $G$.

            \item An ordered graph is a pair $(G, \preccurlyeq)$, where $G$ is a graph and $\preccurlyeq$ is a linear order on $V(G)$. These objects are of interest in extremal combinatorics \cite{orderedgraphs}. An ordered graph can be represented by a structured graph using the map $\sigma \colon V(G)^2 \to 2$, where $\sigma(u, v) = 1$ if and only if $u \preccurlyeq v$.
        \end{enumerate}
    \end{exmps}

    A structured graph $\bm{G} = (G, \sigma)$ is \emphd{Borel} if $G$ is a Borel graph and $\sigma$ is a Borel function (that is, $\sigma^{-1}(i)$ is a Borel subset of $V(G)^{<\infty}$ for all $i \in \N$).

    We naturally extend standard graph-theoretic notation and terminology to structured graphs. For example, given a structured graph $\bm{G} = (G,\sigma)$, we write $V(\bm{G}) \defeq V(G)$ and $E(\bm{G}) \defeq E(G)$. Similarly, given $v \in V(\bm{G})$ and $R \in \N$, we let \[B_{\bm{G}}(v,R) \,\defeq\, B_G(v,R) \,\subseteq\, V(\bm{G}).\] For a subset $U \subseteq V(\bm{G})$, the subgraph of $\bm{G}$ \emphd{induced} by $U$ is the structured graph $\bm{G}[U]$ given by
    \[
        \bm{G}[U] \,\defeq\, \big(G[U], \rest{\sigma}{U^{<\infty}}\big).
    \] 
    Isomorphisms between structured graphs are required to preserve the values of the structure maps. We also let $\asi(\bm{G}) \defeq \asi(G)$.

    With this we can define labelings of structured graphs, rooted (labeled) structured graphs, LCL problems and their solutions, \LOCAL algorithms, and deterministic and randomized \LOCAL complexity of LCL problems on structured graphs by simply repeating Definitions \ref{defn:LCL}, \ref{defn:LOCAL}, \ref{defn:det_comp}, and \ref{defn:rand_comp} with every occurrence of the word ``graph'' replaced by ``structured graph.'' Theorems~\ref{theo:dist_Borel} and \ref{theo:dist_meas} remain true in this more general context \cite[Theorems 2.10 and 2.14]{BernshteynDistributed}, and the same is true for our Theorem~\ref{theo:rand_LOCAL_to_asi}:

    \begin{tcolorbox}
    \begin{theo}[Theorem~\ref{theo:rand_LOCAL_to_asi} for structured graphs]\label{theo:rand_LOCAL_to_asi_structured}
        Let $\Pi$ be an LCL problem and let $\mathbb{G}$ be a class of $\Pi$-colorable finite structured graphs such that $\mathsf{Rand}_{\Pi,\mathbb{G}}(n) = o(\log n)$. Let $\bm{G}$ be a Borel structured graph of finite maximum degree all of whose finite induced subgraphs are in $\mathbb{G}$. If $\asi(\bm{G}) < \infty$, then $\bm{G}$ has a Borel $\Pi$-coloring.
    \end{theo}
    \end{tcolorbox}
    
    In fact, we can strengthen Theorem~\ref{theo:rand_LOCAL_to_asi_structured} by avoiding the asymptotic notation $o(\log n)$. Let $\mathbb{G}$ be a class of finite structured graphs and let $R$, $n \in \N$. A structured graph $\bm{G}$ is \emphd{$(R, n)$-locally in $\mathbb{G}$} if for all $v \in V(\bm{G})$, there is exist an $n$-vertex structured graph $\bm{H} \in \mathbb{G}$ and $u \in V(\bm{H})$ such that
    \[
        \big(\bm{G}[B_{\bm{G}}(v,R)], \, v \big) \,\cong\, \big(\bm{H}[B_{\bm{H}}(u,R)], \, u \big)
    \]
    as rooted structured graphs.

    \begin{tcolorbox}
    \begin{theo}\label{theo:LOCAL_main}
        Let $\Pi = (R, \mathcal{P})$ be an LCL problem and let $\mathbb{G}$ be a class of finite structured graphs. Fix $n \in \N^+$ such that $\mathsf{Rand}_{\Pi, \mathbb{G}}(n) < \infty$ and set $R^\ast \defeq \mathsf{Rand}_{\Pi, \mathbb{G}}(n) + R$.
	Suppose $\bm{G}$ is a Borel structured graph with the following properties:
        \begin{itemize}
            \item $\bm{G}$ is $(R^\ast,n)$-locally in $\mathbb{G}$,

            \item $\asi(\bm{G}) \leq s < \infty$, and

            \item $|B_{\bm{G}}(v, 2R^\ast)| \leq n^{1/(s+1)} / e$ for all $v \in V(\bm{G})$.

        \end{itemize}
        Then $\bm{G}$ has a Borel $\Pi$-coloring.
    \end{theo}
    \end{tcolorbox}
    
    The bound on $|B_{\bm{G}}(v, 2R^\ast)|$ in Theorem~\ref{theo:LOCAL_main} depends on the exact value of $\asi(\bm{G})$. We also prove an analogous statement with a uniform bound for all $\bm{G}$ with $\asi(\bm{G}) < \infty$:  

    \begin{tcolorbox}
    \begin{theo}\label{theo:LOCAL_main_speedup}
        Let $\Pi = (R, \mathcal{P})$ be an LCL problem and let $\mathbb{G}$ be a class of finite  structured graphs. Fix $n \in \N^+$ such that $\mathsf{Rand}_{\Pi, \mathbb{G}}(n) < \infty$ and set $R^\ast \defeq \mathsf{Rand}_{\Pi, \mathbb{G}}(n) + R$.
	Suppose $\bm{G}$ is a Borel structured graph with the following properties:
        \begin{itemize}
            \item $\bm{G}$ is $(R^\ast,n)$-locally in $\mathbb{G}$,

            \item $\asi(\bm{G}) < \infty$, and

            \item $|B_{\bm{G}}(v, 2R^\ast)| \leq n^{1/8} / 4$ for all $v \in V(\bm{G})$.
        \end{itemize}
        Then $\bm{G}$ has a Borel $\Pi$-coloring.
    \end{theo}
    \end{tcolorbox}

    The chain of implications between our main results is as follows:
    \[
        \text{Theorem~\ref{theo:main_easy_asi}} \,\xRightarrow{\text{\ep{\normalfont{A}}}}\,  \text{Theorem~\ref{theo:LOCAL_main}} \,\xRightarrow{\text{\ep{\normalfont{B}}}}\,
        \text{Theorem~\ref{theo:rand_LOCAL_to_asi_structured}} \,\xRightarrow{\text{\ep{\normalfont{C}}}}\, \text{Theorem~\ref{theo:main_speedup}} \,\xRightarrow{\text{\ep{\normalfont{D}}}}\, \text{Theorem~\ref{theo:LOCAL_main_speedup}}.
    \]
    
    Theorem~\ref{theo:main_easy_asi} has been verified in \S\ref{subsec:proof_main_easy}, so now we are left with proving the four implications in the above diagram. In this subsection, we shall establish implications (A), (B), and (D). The only remaining implication, (C), will be proved in \S\ref{subsec:speedup_proof}.

    We start with implication (B), which is proved by a straightforward computation.
    
    \begin{scproof}[ of implication \ep{{\normalfont{B}}}: Theorem~\ref{theo:LOCAL_main} $\Longrightarrow$ Theorem~\ref{theo:rand_LOCAL_to_asi_structured}]\label{proof:B}
        This argument is essentially present in \cite{BernshteynDistributed} (immediately following \cite[Theorem 2.10]{BernshteynDistributed}). Assume that Theorem~\ref{theo:LOCAL_main} holds and let $\Pi = (R, \mathcal{P})$, $\mathbb{G}$, and $\bm{G}$ be as in Theorem~\ref{theo:rand_LOCAL_to_asi_structured}. If $\bm{G}$ is finite, then $\bm{G}$ itself is in $\mathbb{G}$, and so it has an (automatically Borel) $\Pi$-coloring. Thus, we may assume that $\bm{G}$ is infinite. We now claim that the conditions of Theorem~\ref{theo:LOCAL_main} are satisfied for all large enough $n \in \N$, and thus $\bm{G}$ has a Borel $\Pi$-coloring, as desired.
        
        Let the maximum degree of $\bm{G}$ be $\Delta$ and let $s \defeq \asi(\bm{G})$. Let
        \[t(n) \,\defeq\, \mathsf{Rand}_{\Pi, \mathbb{G}}(n) \qquad \text{and} \qquad R^\ast(n) \,\defeq\, t(n) + R.\] Consider any vertex $v \in V(\bm{G})$ and let $B \defeq B_{\bm{G}}(v,R^\ast(n))$. If $n \geq |B|$, then we can extend $B$ to an $n$-element set $U \supseteq B$ of vertices and observe that the $R^\ast(n)$-balls around $v$ in $\bm{G}$ and $\bm{G}[U]$ are the same. Since every finite induced subgraph of $\bm{G}$ is in $\mathbb{G}$, it follows that $G$ is $(R^\ast(n),n)$-locally in $\mathbb{G}$ provided that $n$ is at least as large as the size of every $R^\ast(n)$-ball in $\bm{G}$.
        %
        %
        
        Since $t(n) = o(\log n)$ and hence $\Delta^{t(n)} = n^{o(1)}$, the following chain of inequalities holds for all sufficiently large $n \in \N$: 
    \begin{equation}\label{eq:bounds_for_n} 
        1 + \Delta^{R^*(n)} \,\leq\, 1+\Delta^{2R^*(n)} \,=\, 1 + \Delta^{2t(n) + 2R} \,\leq\, n^{1/(s+1)}/e \,\leq\, n.
    \end{equation}
    As every $r$-ball in $\bm{G}$ contains at most $1 + \Delta^r$ vertices, any $n$ satisfying the bounds in \eqref{eq:bounds_for_n} fulfills all the requirements 
    of Theorem~\ref{theo:LOCAL_main}, and we are done.
    \end{scproof}

    To prove implications (A) and (D), we rely on the following technical lemma from \cite{BernshteynDistributed}: 

    \begin{lemma}[{\cite[Lemma 4.6]{BernshteynDistributed}}]\label{lemma:LCL_to_CSP}
    Let $\Pi = (R, \mathcal{P})$ be an LCL problem and let $\mathbb{G}$ be a class of finite structured graphs. Fix $n \in \N^+$ such that $T \defeq \mathsf{Rand}_{\Pi, \mathbb{G}}(n)$ is finite and set $R^\ast \defeq T + R$. Let $\ell \in \N^+$ and a \LOCAL algorithm $\mathcal{A}$ witness the bound $\mathsf{Rand}_{\Pi, \mathbb{G}}(n) \leq T$ \ep{as in Definition~\ref{defn:rand_comp}}. If $\bm{G}$ is a Borel structured graph of finite maximum degree that is $(R^\ast, n)$-locally in $\mathbb{G}$, then there exists a bounded Borel CSP $\B \colon V(\bm{G}) \to^? \ell$ such that: 
		\begin{enumerate}[label=\ep{\normalfont\roman*}]
			\item for every solution $\theta \colon V(\bm{G}) \to \ell$ to $\B$, the function $\mathcal{A}_T(\bm{G}, \theta)$ is a $\Pi$-coloring of $\bm{G}$, and
			
			\item $\pr(\B) \leq 1/n$ and $\de(\B) \leq \sup_{v \in V(\bm{G})} |B_{\bm{G}}(v, 2R^\ast)| - 1$.
		\end{enumerate}
        Moreover, the associated graph $G_\B$ is a subgraph of $G^{2R^\ast}$, and hence $\asi(G_\B) \leq \asi(\bm{G})$.
    \end{lemma}

    We point out that the ``moreover'' part of Lemma~\ref{lemma:LCL_to_CSP} is not explicitly stated in \cite{BernshteynDistributed}, but it is an immediate consequence of the proof of \cite[Lemma 4.6]{BernshteynDistributed}, since the domain of each constraint in the CSP $\B$ constructed there is the $R^\ast$-ball around some vertex of $\bm{G}$.

    \begin{scproof}[ of implication \ep{{\normalfont{A}}}: Theorem~\ref{theo:main_easy_asi} $\Longrightarrow$ Theorem~\ref{theo:LOCAL_main}]
        Let $\Pi = (R, \mathcal{P})$, $\mathbb{G}$, $n \in \N^+$, $R^\ast$, $\bm{G}$, and $s$ be as in Theorem~\ref{theo:LOCAL_main}. Set $T \defeq \mathsf{Rand}_{\Pi, \mathbb{G}}(n)$, so that $R^\ast = T + R$. Let $\ell \in \N^+$ and a \LOCAL algorithm $\mathcal{A}$ witness the bound $\mathsf{Rand}_{\Pi, \mathbb{G}}(n) \leq T$. Since $\bm{G}$ is $(R^\ast, n)$-locally in $\mathbb{G}$, Lemma~\ref{lemma:LCL_to_CSP} yields a bounded Borel CSP $\B \colon V(\bm{G}) \to^? \ell$ such that:
        \begin{enumerate}[label=\ep{\normalfont\roman*}]
			\item\label{item:theta} for every solution $\theta \colon V(\bm{G}) \to \ell$ to $\B$, the function $\mathcal{A}_T(\bm{G}, \theta)$ is a $\Pi$-coloring of $\bm{G}$,
			
			\item\label{item:pdbounds} $\pr(\B) \leq 1/n$ and $\de(\B) \leq \sup_{v \in V(\bm{G})} |B_{\bm{G}}(v, 2R^\ast)| - 1  \leq n^{1/(s+1)}/e - 1$, and

            \item\label{item:asibound} $\asi(G_\B) \leq s$.
		\end{enumerate}
        It follows from \ref{item:pdbounds} that
        \[
            \pr(\B) \, (\de(\B) + 1)^{s+1} \,\leq\, \frac{1}{n} \,\cdot\, \left(n^{1/(s+1)}/e\right)^{s+1} \,=\, e^{-s - 1}.
        \]
        Therefore, \ref{item:asibound} implies that $\B$ has a Borel solution $\theta \colon V(\bm{G}) \to \ell$ by Theorem~\ref{theo:main_easy_asi}. Then, by \ref{item:theta},
        $
            f \defeq \mathcal{A}_T(\bm{G}, \theta)
        $
        is a $\Pi$-coloring of $\bm{G}$, and it is Borel by \cite[remark after Lemma 4.1]{BernshteynDistributed}.
    \end{scproof}

    Implication (D), i.e., Theorem~\ref{theo:main_speedup} $\Longrightarrow$ Theorem~\ref{theo:LOCAL_main_speedup}, is proved in exactly the same way, \emph{mutatis mutandis}.  

    \subsection{From \LOCAL algorithms to the LLL}\label{subsec:speedup_proof}

    In this section we verify the only remaining implication from the diagram in \S\ref{subsec:LOCAL_proof}, namely (C): Theorem~\ref{theo:rand_LOCAL_to_asi_structured} $\Longrightarrow$ Theorem~\ref{theo:main_speedup}. Our argument follows the strategy outlined in \cite[\S5.D]{BernshteynDistributed}. The key point is that there exist efficient randomized \LOCAL algorithms that solve CSPs satisfying a \pw LLL criterion, and combining such an algorithm with Theorem~\ref{theo:rand_LOCAL_to_asi_structured} yields Borel solutions to such CSPs when the associated graph has finite asymptotic separation index.

    To make the above idea precise, we need to be able to encode a given CSP as an LCL problem on an auxiliary structured graph. There are several natural ways of doing so; here we follow the approach described in \cite[\S5.D]{BernshteynDistributed}.
    
    \begin{defn}[Graph-CSPs]
        A \emphd{graph-CSP} with range $q \in \N^+$ is a pair $(G, \B)$, where $G$ is a graph and $\B \colon V(G) \to^? q$ is a bounded CSP such that $\de(\B) < \infty$ and $G_\B$ is a subgraph of $G$.
    \end{defn}
    
    A graph-CSP $(G, \B)$ can be viewed as a structured graph in the following way. For each tuple of vertices $\bm{v} = (v_0, \ldots, v_{k-1}) \in V(G)^{<\infty}$, the structure map $\sigma$ will contain the information about the constraints $B \in \B$ such that $\dom(B) = \set{v_0, \ldots, v_{k-1}}$. Formally, let $(V(G))^{< \infty}$ be the set of all finite tuples of distinct vertices and for each $\bm{v} = (v_0, \ldots, v_{k-1}) \in (V(G))^{< \infty}$, define
    \[
        \B_{\bm{v}} \,\defeq\, \set{B \in \B \,:\, \dom(B) = \set{v_0, \ldots, v_{k-1}}}.
    \]
    The domain of $\sigma$ is the set of all $\bm{v} \in (V(G))^{<\infty}$ with $\B_{\bm{v}} \neq \0$. Note that since $\B$ is bounded, $\sigma$ is defined on tuples of bounded length. For $\bm{v} = (v_0, \ldots, v_{k-1}) \in \dom(\sigma)$, let $\iota_{\bm{v}} \colon \set{v_0, \ldots, v_{k-1}} \to k$ be the bijection given by $\iota_{\bm{v}}(x_i) \defeq i$, and for each $B \in \B_{\bm{v}}$, define
    \[
        B_{\bm{v}}^\ast \, \defeq\,  \set{ \phi \colon k \to q \,:\, \phi \circ \iota_{\bm{v}} \in B}.
    \]
    Note that the mapping $B^\ast_{\bm{v}} \to B \colon \phi \mapsto \phi \circ \iota_{\bm{v}}$ is a bijection. Now let the \emphd{type} of $\bm{v}$ be
 	\[
 		\mathsf{type}(\bm{v}) \,\defeq\, \set{B^\ast_{\bm{v}} \,:\, B \in \B_{\bm{v}}}.
 	\]
    Observe that the set $\mathsf{Types}$ of all possible types of tuples in $\dom(\sigma)$ is countable (indeed, it is a subset of $\finset{\finset{\finfun{\N}{\N}}}$), so we can fix an arbitrary injection $\mathsf{code} \colon \mathsf{Types} \to \N$ and let
  \[
    \sigma(\bm{v}) \,\defeq\, \mathsf{code}(\mathsf{type}(\bm{v})).
    \]
    The pair $(G, \sigma)$ now represents exactly the same information as $(G, \B)$, and since $\sigma$ is defined on tuples of bounded length, $(G, \sigma)$ is a structured graph (as in Definition~\ref{defn:str_graph}). Moreover, it is routine to check that if $G$ is a Borel graph and $\B$ is a Borel CSP, then $\sigma$ is a Borel function, and thus $(G,\sigma)$ is a Borel structured graph.

    Let $(G,\mathscr{B})$ be a graph-CSP. The requirement that $G_\B$ is a subgraph of $G$ means that the $1$-ball around every vertex of $G$ includes the domains of all constraints that involve that vertex. That is, the problem of solving a graph-CSP is an LCL problem of radius $1$, which we denote by $\Pi_{\mathrm{CSP}, q}$ (here $q \in \N^+$ is the range of $\B$). Given parameters $q$, $k$, $p$, and $d$, let $\mathbb{CSP}(q,k,p,d)$ be the class of all finite graph-CSPs $(G, \B)$ (viewed as structured graphs) such that:
 	\begin{itemize}
 		\item the range of $\B$ is $q$;
 		\item $\sup_{B \in \B} |\dom(B)| \leq k$;
 		\item $\pr(\B) \leq p$ and $\de(\B) \leq d$.
 	\end{itemize}
  
 	After these technical preliminaries, we can rigorously formulate a result regarding \LOCAL algorithms for the LLL that we shall use: 
 	
 	\begin{theo}[{Ghaffari--Harris--Kuhn \cite{GHK}}]\label{theo:dist_LLL}
 		If $q$, $k$, $p$, and $d$ satisfy $p(d+1)^8 \leq 2^{-15}$, then
 		\[
 			\mathsf{Rand}_{\Pi_{\mathrm{CSP}, q},\, \mathbb{CSP}(q,k,p,d)}(n) \,=\, \exp(O(\sqrt{\log \log n})) \,=\, o(\log n).
 		\]
 	\end{theo}

        The first $o(\log n)$-time randomized \LOCAL algorithm for the LLL with a \pw criterion is due to Fischer and Ghaffari \cite{FG}. We are using the Ghaffari--Harris--Kuhn version instead because it improves the criterion from $\pr(\de+1)^{32} \leq e^{-32}$ to $\pr(\de+1)^8 \leq 2^{-15}$. In subsequent work, considerable effort has been put into making the running time of the algorithm even lower \cite{RG,Davies}. On the other hand, we are not aware of any improvements to the required bound on $\pr$ as a function of $\de$, and it remains an open problem whether the exponent $8$ can be reduced, ideally all the way down to $1$.

        With Theorem~\ref{theo:dist_LLL} in hand, we are ready to deduce Theorem~\ref{theo:main_speedup} from Theorem~\ref{theo:rand_LOCAL_to_asi_structured}:

        \begin{scproof}[ of implication \ep{{\normalfont{C}}}: Theorem~\ref{theo:rand_LOCAL_to_asi_structured} $\Longrightarrow$ Theorem~\ref{theo:main_speedup}]
            Let $\B \colon X \to^? q$ be a bounded Borel CSP on a standard Borel space $X$ such that $\asi(G_\B) < \infty$ and
            \[
                \pr(\B) \, (\de(\B) + 1)^8 \,\leq\, 2^{-15}.
            \]
        Let $p \defeq \pr(\B)$, $d \defeq \de(\B)$, and $k \defeq \sup_{B \in \B} |\dom(B)|$, and let $\Pi \defeq \Pi_{\mathrm{CSP}, q}$ and $\mathbb{G} \defeq \mathbb{CSP}(q,k,p,d)$. By the usual LLL, the structured graphs in $\mathbb{G}$ are $\Pi$-colorable. Furthermore, by Theorem~\ref{theo:dist_LLL}, we have $\mathsf{Rand}_{\Pi, \mathbb{G}}(n) = o(\log n)$.   Consider the graph-CSP $\bm{G} \defeq (G_\B, \B)$, viewed as a Borel structured graph. The maximum degree of $\bm{G}$ is bounded above by $kd$, and all finite induced subgraphs of $\bm{G}$ are in $\mathbb{G}$ (taking an induced subgraph corresponds to passing to a subset of the vertices and only retaining the constraints whose domains are contained in that subset). Since $\asi(\bm{G}) < \infty$ by assumption, we conclude that $\bm{G}$ has a Borel $\Pi$-coloring---i.e., $\B$ has a Borel solution---by Theorem~\ref{theo:rand_LOCAL_to_asi_structured}.
        \end{scproof}

    \section{Applications}

    \subsection{Coloring graphs without short cycles}\label{subsec:Joh_proof}

    Here we present the derivation of Theorem~\ref{theo:BorelJohansson} from our main results.  For $\Delta \in \N$, let $\mathbb{T}\mathrm{r}\mathbb{F}\mathrm{ree}_\Delta$ be the class of all triangle-free finite graphs of maximum degree at most $\Delta$. Given $\epsilon > 0$ and $\Delta \geq 2$, let
    \[
        q(\Delta, \epsilon) \,\defeq\, \left\lfloor (4 + \epsilon) \frac{\Delta}{\log\Delta} \right\rfloor.
    \]
    Let $\Pi_q$ denote the LCL problem that encodes proper $q$-coloring of graphs (see Example~\ref{exmp:coloringLCL}). Chung, Pettie, and Su \cite[Theorem 9]{CPS} showed that for each $\epsilon > 0$, there exist $C$, $\delta$, $\Delta_0 > 0$ (depending on $\epsilon$) such that for every $\Delta \geq \Delta_0$, we have
    \begin{equation}\label{eq:CPS}
        \mathsf{Rand}_{\Pi_{q(\Delta, \epsilon)}, \, \mathbb{T}\mathrm{r}\mathbb{F}\mathrm{ree}_\Delta}(n) \,\leq\, C\Delta^{-\delta} \log n.
    \end{equation}
    Although this bound is not strong enough to invoke Theorem~\ref{theo:rand_LOCAL_to_asi} (it is of order $O(\log n)$ rather than $o(\log n)$), it is still sufficient for a direct application of Theorem~\ref{theo:LOCAL_main_speedup}.\footnote{Actually, the running time of the Chung--Pettie--Su algorithm can be improved to $o(\log n)$ via generic ``speed-up'' theorems in distributed computing, namely by combining a result of Chang and Pettie \cite[Theorem 6]{CP} with an $o(\log n)$-time randomized \LOCAL algorithm for the LLL such as Theorem~\ref{theo:dist_LLL}. That being said, the bound \eqref{eq:CPS} suffices for our purposes.} Indeed, the same analysis as in the proof of the implication Theorem~\ref{theo:LOCAL_main} $\Longrightarrow$ Theorem~\ref{theo:rand_LOCAL_to_asi_structured} on page \pageref{proof:B} reduces the problem to showing that for all sufficiently large $n \in \N$,
    \[
        1 + \Delta^{2C\Delta^{-\delta} \log n + 2} \,\leq\, \frac{1}{4}n^{1/8}. 
    \]
    This is indeed the case as long as $\Delta$ is large enough as a function of $\epsilon$, since
    \[
        \Delta^{2C\Delta^{-\delta} \log n} \,=\, n^{2C\Delta^{-\delta} \log \Delta},
    \]
    and $2C\Delta^{-\delta} \log \Delta < 1/8$ for sufficiently large $\Delta$. This establishes the part of Theorem~\ref{theo:BorelJohansson} regarding triangle-free graphs. The part of Theorem~\ref{theo:BorelJohansson} concerning graphs without $3$- and $4$-cycles is proved in exactly the same way, but with \cite[Theorem 9]{CPS} replaced by the randomized \LOCAL algorithm of Chung, Pettie, and Su for $\lfloor(1+\epsilon)\Delta / \log \Delta\rfloor$-coloring graphs without cycles of length at most $4$ \ep{see the remark in \cite{CPS} after \cite[Theorem 9]{CPS}}.

    \subsection{Edge-colorings with $\Delta + \tilde{O}(\sqrt{\Delta})$ colors}\label{subsec:proof_edge_epsilon}

    In this subsection, we prove Theorem~\ref{theo:Borel_edge-coloring_close_to_Delta}. We begin by explaining how to frame proper edge-coloring as an LCL problem. The obvious approach is to view proper $q$-edge-coloring of a graph $G$ as a proper $q$-coloring of the line graph $L(G)$ and use the LCL problem defined in Example~\ref{exmp:coloringLCL}. The drawback of this strategy is that \LOCAL algorithms for edge-coloring typically reference both edges and vertices of $G$, and the information about the vertices of $G$ is apparently lost when working on $L(G)$. However, this issue can be easily rectified by incorporating both the vertices and the edges of $G$ into an auxiliary structure. For example, we can use the following construction, suggested in \cite[5]{ChangReduction}. Let $G$ be a graph. Define a new graph $G'$ with vertex set $V(G') \defeq V(G) \sqcup E(G)$ and edge set $E(G') \defeq \set{\set{v,e} \,:\, v \in V(G),\, e \in E(G), \, v \in e}$. For all $v \in V(G)$ and $e \in E(G)$, let $\sigma(v) \defeq 0$ and $\sigma(e) \defeq 1$. This defines a function $\sigma \colon V(G') \to 2$ and yields a structured graph $\bm{G'} \defeq (G', \sigma)$, which contains information about the vertices and edges of $G$, distinguished from each other via the map $\sigma$. Proper $q$-edge-coloring of $G$ can be encoded as an LCL problem on $\bm{G'}$ of radius $2$, which requires every vertex $x$ of $\bm{G'}$ with $\sigma(x) = 1$ to receive a label in $q$ that is distinct from the labels of the vertices joined to $x$ by paths of length $2$.

    Fix $\Delta \in \N$ and let 
    $
        q \defeq \lfloor \Delta + \sqrt{\Delta} \log^3\Delta\rfloor
    $.
    Combining a result of Chang, He, Li, Pettie, and Uitto \cite[Theorem 4]{CHLPU} with a randomized LLL algorithm due to Fischer and Ghaffari \cite{FG} (or with Theorem~\ref{theo:dist_LLL}), we see that the proper $q$-edge-coloring problem can be solved on $n$-vertex graphs of maximum degree $\Delta$ by a randomized \LOCAL algorithm with running time $o(\log n)$, assuming $\Delta$ is sufficiently large. 
    By Theorem~\ref{theo:LOCAL_main_speedup}, it follows that every Borel graph $G$ with finite asymptotic separation index and of maximum degree $\Delta$ has a Borel proper $q$-edge-coloring, as desired.
    
    For further discussion of the state of the art in distributed edge-coloring, see \cite{Davies}.


    \subsection{Borel independent complete sections}\label{subsec:proof_cs}

    Here we prove Theorem~\ref{theo:indep_transversal}, restated below for ease of reference.

    \begin{theocopy}{theo:indep_transversal}\label{theo:indep_transversal1}
        Let $H$ be a Borel graph with $\asi(H) < \infty$ and let $G_1$, $G_2$ be two spanning Borel subgraphs of $H$. Suppose that the maximum degree of $G_1$ is at most $\Delta \in \N$, where $\Delta \geq 2$. 
        If every component of $G_2$ contains at least $
            2^{25} \,\Delta\log\Delta
        $ 
        vertices, then there is a Borel $G_1$-independent complete section for $G_2$.
    \end{theocopy}
    \begin{scproof}[ of Theorem~\ref{theo:indep_transversal1}]
        For brevity, let $V \defeq V(H)$. Let $k \defeq \lceil 2^{25} \Delta \log\Delta \rceil$, so every component of $G_2$ has at least $k$ vertices. 
        %
        In fact, we may assume that every component of $G_2$ contains \emph{exactly} $k$ vertices. Indeed, let $\mathcal{R} \subseteq \finset{V}$ be the family of all finite sets $R$ such that:
        \begin{itemize}
            \item $|R| = k$, and
            \item the graph $G_2[R]$ is connected.
        \end{itemize}
        By \cite[Lemma 7.3]{KechrisMiller}, there exists a Borel set $\mathcal{R}' \subseteq \mathcal{R}$ such that the members of $\mathcal{R}'$ are pairwise disjoint and every set in $\mathcal{R}$ intersects some member of $\mathcal{R}'$. It follows that the vertex set of every connected component of $G_2$ includes at least one member of $\mathcal{R}'$ as a subset. Now define $U \defeq \bigcup \mathcal{R}'$, $H' \defeq H[U]$, $G_1' \defeq G_1[U]$, and let $G_2'$ be the spanning subgraph of $G_2[U]$ obtained by only retaining the edges of the graphs $G_2[R]$ with $R \in \mathcal{R}'$ (in other words, $G_2'$ is the graph whose connected components are precisely the graphs $G_2[R]$ with $R \in \mathcal{R}'$). Then $H'$, $G_1'$, and $G_2'$ are Borel graphs satisfying the assumptions of Theorem~\ref{theo:indep_transversal1}. Furthermore, by construction, every component of $G_2'$ contains exactly $k$ vertices. Suppose we can find a Borel $G_1'$-independent complete section $S$ for $G_2'$. Then $S$ is $G_1$-independent (because $G_1'$ is an induced subgraph of $G_1$) and a complete section for $G_2$ (because every component of $G_2$ includes a component of $G_2'$). Therefore, proving the theorem under the assumption that every component of $G_2$ contains exactly $k$ vertices implies the theorem in general.

        Let $G \defeq G_1$ and $\mathcal{F} \defeq \set{V(C) \,:\, \text{$C$ is a connected component of $G_2$}}$. Then $\mathcal{F}$ is a Borel subset of $[V]^k$, and our goal is to find a Borel $G$-independent complete section for $\mathcal{F}$. Without loss of generality, we may assume that each set $F \in \mathcal{F}$ is $G$-independent. To see this, let $G'$ be the spanning subgraph of $G$ in which two vertices are adjacent if and only if they are adjacent in $G$ and do not belong to the same set in $\mathcal{F}$ (so each set $F \in \mathcal{F}$ is $G'$-independent by construction). Suppose we can find a Borel $G'$-independent complete section $S$ for $\mathcal{F}$. Using the Luzin--Novikov theorem \ep{Theorem~\ref{theo:LN}}, we can choose, for each $F \in \mathcal{F}$, a single vertex $v_F \in S \cap F$ so that the mapping $\mathcal{F} \to S \colon F \mapsto v_F$ is Borel. Then $S' \defeq \set{v_F \,:\, F \in \mathcal{F}}$ is a Borel subset of $S$ that meets every set in $\mathcal{F}$ in exactly one vertex. It follows that $S'$ is $G$-independent, so we have found a Borel $G$-independent complete section for $\mathcal{F}$.

        With these preparations, we may start the main part of the proof. Given a $\Delta$-coloring $f \colon V \to \Delta$, we define a set $S_f \subseteq V$ by 
        \[
            S_f \,\defeq\, \set{v \in V \,:\, f(v) = 0 \text{ and } f(u) \neq 0 \text{ for all } uv \in E(G)}.
        \]
        By construction, $S_f$ is a $G$-independent set, and if $f$ is a Borel function, then $S_f$ is Borel. Our aim is to find a Borel function $f \colon V \to \Delta$ such that $S_f$ meets every set in $\mathcal{F}$.

        The membership of a vertex $v$ in $S_f$ is determined by the values of $f$ at $v$ and the neighbors of $v$. Thus, whether $S_f$ meets a set $F \in \mathcal{F}$ only depends on the restriction of $f$ to $B_G(F, 1)$, i.e., the set of all vertices that are in $F$ or have a $G$-neighbor in $F$. Hence, we can define a $(V,\Delta)$-constraint $B_F$ with domain $B_G(F,1)$ such that $f \colon V \to \Delta$ satisfies $B_F$ if and only if $S_f$ meets $F$. Explicitly, 
        \[
            B_F \,\defeq\, \big\{\phi \colon B_G(F,1) \to \Delta \ : \ \text{for all $v \in F$, $\phi(v) \neq 0$ or $\phi(u) = 0$ for some $uv \in E(G)$}\big\}.
        \]
        Let $\B \defeq \set{B_F \,:\, F \in \mathcal{F}}$. This is a bounded Borel CSP such that if $f \colon V \to \Delta$ is a solution to $\B$, then $S_f$ is a Borel $G$-independent complete section for $\mathcal{F}$. It remains to argue that $\B$ has a Borel solution, to which end we shall employ Theorem~\ref{theo:main_speedup}.

        To begin with, we observe that $\asi(G_\B) < \infty$. Indeed, if $u$ and $v$ are neighbors in $G_\B$, then there is some $F \in \mathcal{F}$ such that $u$, $v \in B_G(F,1)$. Since $H[F]$ is a connected graph with at most $k$ vertices, it follows that $u$ and $v$ are joined by a path of at most $k+1$ edges in $H$. Thus, $G_\B$ is a subgraph of $H^{k+1}$, and hence $\asi(G_\B) \leq \asi(H) < \infty$. We will not use the graph $H$ again, so all the references to neighbors, adjacency, edges, etc.~in the sequel are in relation to the graph $G$.

        Next we need to bound the parameters $\de(\B)$ and $\pr(\B)$. We start with $\de(\B)$, which is the easier one. Take $F \in \mathcal{F}$ and suppose that $F' \in \mathcal{F} \setminus \set{F}$ satisfies $B_G(F',1) \cap B_G(F,1) \neq \0$. Consider any vertex $v \in B_G(F,1) \cap B_G(F',1)$. Since $v \in B_G(F,1)$, there are at most $k + \Delta k$ choices for $v$. As $v \in B_G(F',1)$, either $v$ itself or a neighbor of $v$ is in $F'$. Since the sets in $\mathcal{F}$ are disjoint, it follows that, given $v$, there are at most $1 + \Delta$ choices for $F'$. In total, for fixed $F$, there are at most $(k + \Delta k)(1 + \Delta) - 1$ options for $F'$ (we are subtracting $1$ since $F$ itself should not be counted), so
        \[
            \de(\B) \, \leq \, (k + \Delta k)(1 + \Delta) - 1 \,\leq\, 4\Delta^2 k - 1 \,\leq\, 2^{27} \Delta^5 - 1.
        \]

        Now we turn our attention to $\pr(\B)$. Consider any set $F \in \mathcal{F}$. We need to bound $\P[B_F]$, i.e., the probability that a uniformly random $\Delta$-coloring $\phi \colon B_G(F,1) \to \Delta$ belongs to the set $B_F$. Let $N$ be the set of all neighbors of the vertices in $F$, and let $\psi \colon N \to \Delta$ and $\xi \colon F \to \Delta$ be the restrictions of $\phi$ to $N$ and $F$ respectively. By our assumption, the set $F$ is $G$-independent, so, $B_G(F,1) = F \sqcup N$, which means that $\psi$ and $\xi$ are chosen independently and uniformly at random from the sets of all $\Delta$-colorings of $N$ and $F$ respectively. Let
        \begin{equation}\label{eq:Sast_defn}
            F^\ast \,\defeq\, \set{v \in F \,:\, \psi(u) \neq 0 \text{ for all neighbors } u \text{ of }v}.
        \end{equation}
        Note that $\phi \in B_F$ if and only if $\xi(v) \neq 0$ for all $v \in F^\ast$. Therefore,
        \begin{equation}\label{eq:BS_bound}
            \P[B_F] \,=\, \P[\xi(v) \neq 0 \text{ for all } v \in F^\ast].
        \end{equation}
        The key observation is that $|F^\ast|$ is unlikely to be too small:

        \begin{claim}\label{claim:Qast}
            $
             \displaystyle  \P[|F^\ast| < k/8] \,\leq\, 4\Delta^{-500}. 
            $
        \end{claim}
        \begin{claimproof}
            Observe that, since $(1-1/x)^x \geq (1-1/2)^2 = 1/4$ for all $x \geq 2$, we have
            \[
                \E[|F^\ast|] \,\geq\, |F| \left(1 - \frac{1}{\Delta}\right)^\Delta \,\geq\, \frac{k}{4}.
            \]
            Therefore, we need to find an upper bound on the probability that $|F^\ast| < \E[|F^\ast|] - k/8$. To this end, we shall apply a concentration of measure result that is often used in probabilistic combinatorics and is a consequence of Talagrand's isoperimetric inequality \cite{Talagrand}. There are several variants of this inequality in the literature \cite{MolloyReedDelta, BruhnJoos, KellyPostle}, and designing new and improved versions of it is an active area of research. We shall use a (simplified form of) the recent formulation due to Delcourt and Postle \cite{DelcourtPostle}. Given two finite strings ${x} = (x_1, \ldots, x_n)$ and ${y} = (y_1, \ldots, y_n)$, and a vector ${c} = (c_1, \ldots, c_n) \in \R^n$, we let
            \[
                \dist_{{c}}({x}, {y}) \,\defeq\, \sum_{i \,:\, x_i \neq y_i} |c_i|.
            \]

        \begin{theo}[{Talagrand's inequality; Delcourt--Postle \cite[Theorem 4.4]{DelcourtPostle}}]\label{theo:Talagrand}
            Let $(\Omega_1, \mu_1)$, \ldots, $(\Omega_n, \mu_n)$ be discrete probability spaces and let $(\Omega, \mu) \defeq \prod_{i = 1}^n (\Omega_i, \mu_i)$ be their product space. Suppose that $X \colon \Omega \to \R$ is a non-negative random variable satisfying the following for some $r$, $d > 0$:

            \begin{enumerate}[label=\ep{\normalfont{T}}]
                \item\label{item:T} For each 
                ${\omega} \in \Omega$, 
                there exists a vector ${c} \in \R^n$ such that $\|{c}\|_1 \leq rX(\omega)$, $\|{c}\|_\infty \leq d$, and for all ${\omega}' \in \Omega$, we have $X({\omega}') \geq X(\omega) - \dist_{{c}}({\omega}',{\omega})$.
            \end{enumerate}
            Then for any $t > 96 \sqrt{rd\E[X]} + 128rd$,
            \[
                \P\left[|X - \E[X]| > t \right] \,\leq\, 4 \exp\left(- \frac{t^2}{8rd(4\E[X] + t)}\right).
            \]
        \end{theo}

        To apply this in our setting, we view $\psi \colon N \to \Delta$ as a point in the product of $|N|$ copies of the space $(\Delta, \nu)$, where $\nu$ is the uniform probability measure on the set $\Delta = \set{0,1,\ldots, \Delta - 1}$. Let \[X \,\defeq\, |F| - |F^\ast|,\] 
        where $F^*(\psi)$ is constructed according to \eqref{eq:Sast_defn} using $\psi$.
        We claim that condition \ref{item:T} of Theorem~\ref{theo:Talagrand} holds for the random variable $X$ with $r = 1$ and $d = \Delta$. 
        For each vertex $v \in F \setminus F^\ast(\psi)$, we arbitrarily pick a neighbor $u_v \in N$ such that $\psi(u_v) = 0$ (which exists by the definition of $F^\ast(\psi)$) 
        and define a function $c \colon N \to \N$ as follows:
        \[
            c(u) \,\defeq\, |\set{v \in F \setminus F^\ast(\psi) \,:\, u = u_v}|.
        \]
        We have $\|c\|_1 = \sum_{u \in N} c(u) = X(\psi)$ and $\|c\|_\infty \leq \Delta$, since each $u \in N$ has at least $c(u)$ neighbors in $F \setminus F^\ast(\psi)$. 
        Now take any $\psi' \colon N \to \Delta$ and note that if $v \in F^\ast(\psi') \setminus F^\ast(\psi)$, then $\psi'(u_v) \neq 0$, so
        \[
            X(\psi') \,\geq\,  X(\psi) - |F^\ast(\psi') \setminus F^\ast(\psi)| \,\geq\, X(\psi) - \sum_{u \,:\, \psi'(u) \neq 0} c(u) \,=\, X(\psi') - \dist_c(\psi, \psi'),
        \]
        as desired. 
            Since $\E[X] \leq k$ and $k/8 > 96\sqrt{\Delta k} + 128\Delta$, we may now apply Theorem~\ref{theo:Talagrand} with $r=1$, $d=\Delta$, and $t = k/8$ to conclude that
            \begin{align*}
                \P\big[|F^\ast| < \E[|F^\ast|] -  k/8\big] \,&=\, \P[X - \E[X] > k/8] \\
                &\leq\, 4 \exp\left(-\frac{k^2/64}{8\Delta(4k + k/8)}\right) \,=\, 4\exp\left( - \frac{k}{2112\Delta}\right) \,\leq\, 4\Delta^{-500}, 
            \end{align*}
            where the last inequality uses the bound $k \geq 2^{25}\,\Delta \log\Delta$.
        \end{claimproof}

        Note that the definition of the set $F^\ast$ depends on $\psi$ but not on $\xi$. Since the functions $\psi \colon N \to \Delta$ and $\xi \colon F \to \Delta$ are chosen independently, conditioned on the choice of $\psi$, we can write
        \[
            \P[\xi(v) \neq 0 \text{ for all } v \in F^\ast \,\vert\, \psi] \,=\, \left(1 - \frac{1}{\Delta}\right)^{|F^\ast|}.
        \]
        Therefore, thanks to Claim~\ref{claim:Qast}, we have
        \begin{align*}
            \P[\xi(v) \neq 0 \text{ for all } v \in F^\ast] \,&\leq\, 4\Delta^{-500} \,+\, \left(1 - \frac{1}{\Delta}\right)^{k/8}\\
            \big(\text{since } 1-x \leq e^{-x} \text{ for all } x \in \R\big) \qquad \qquad &\leq\, 4\Delta^{-500}\,+\, \exp\left(-\frac{k}{8\Delta}\right) \,\leq\, 5\Delta^{-500},
        \end{align*}
        where the last inequality uses the bound $k \geq 2^{25}\Delta \log\Delta$. Remembering \eqref{eq:BS_bound}, we conclude that
        \[
            \P[B_F] \,\leq\, 5\Delta^{-500}.
        \]
        Since this bound holds for all $F \in \mathcal{F}$, it follows that $\pr(\B) \leq 5\Delta^{-500}$.

        Putting everything together, we see that
        \[
            \pr(\B) \, (\de(\B) + 1)^8 \,\leq\, 5 \Delta^{-500} \cdot 2^{216} \Delta^{40} \,=\, 5 \cdot 2^{216} \cdot \Delta^{-460} \,\leq\, 5 \cdot 2^{-244} \,<\, 2^{-15},
        \]
        and thus $\B$ has a Borel solution by Theorem~\ref{theo:main_speedup}, as desired.
    \end{scproof}

    \subsection{Edge-colorings of Schreier graphs}\label{subsec:proof_Schreier}

    In this subsection, we derive Theorem~\ref{theo:BorelVizing} from Theorem~\ref{theo:indep_transversal}. The main idea of the proof of Theorem~\ref{theo:BorelVizing} is to decompose the line graph $L(G)$ of $G$ into certain induced subgraphs with a simple structure and then combine colorings of those subgraphs into a coloring of $L(G)$. We begin with a couple lemmas that facilitate this approach. The following observation is almost trivial:
    
    \begin{lemma}\label{lemma:decomposition_simple}
        Let $G$ be a Borel graph and let $V(G) = U_1 \sqcup \ldots \sqcup U_n$ be a partition of the vertex set of $G$ into Borel subsets such that $\chi_\mathsf{B}(G[U_i]) < \infty$ for all $1 \leq i \leq n$. Then
        \[
            \chi_\mathsf{B}(G) \,\leq\, \sum_{i=1}^n \chi_\mathsf{B}(G[U_i]).
        \]
    \end{lemma}
    \begin{scproof}
        Let $q_i \defeq \chi_\mathsf{B}(G[U_i])$ and let $f_i \colon U_i \to q_i$ be a Borel proper $q_i$-coloring of $G[U_i]$. Define
        \[
            f_i'(v) \,\defeq\, f_i(v) + \sum_{j=1}^{i-1} q_j \quad \text{for all } v \in U_i.
        \]
        Then $f \defeq f_1' \sqcup \ldots \sqcup f_n'$ is a Borel proper $\left(\sum_{i=1}^n q_i\right)$-coloring of $G$.
    \end{scproof}

    For a graph $G$, let $\Delta(G)$ denote the maximum degree of $G$. If in the setting of Lemma~\ref{lemma:decomposition_simple} each graph $G[U_i]$ has finite maximum degree, then we can use Theorem~\ref{theo:KST} to write
    \[
        \chi_\mathsf{B}(G) \,\leq\, \sum_{i=1}^n \chi_\mathsf{B}(G[U_i]) \,\leq\, \sum_{i=1}^n(\Delta(G[U_i]) + 1) \,=\, n + \sum_{i=1}^n\Delta(G[U_i]). 
    \]
    The next lemma, which we deduce from Theorem~\ref{theo:indep_transversal}, allows us to improve this bound when $G$ has finite asymptotic separation index and the graphs $G[U_i]$ do not have small components:
    
    \begin{lemma}\label{lemma:decomposition_complicated}
        Let $G$ be a Borel graph of maximum degree at most $\Delta \in \N$, where $\Delta \geq 2$, and let $V(G) = U_1 \sqcup \ldots \sqcup U_n$ be a partition of the vertex set of $G$ into Borel subsets. Suppose that for all $1 \leq i \leq n$, every component of $G[U_i]$ contains at least $2^{25} \,\Delta\log\Delta$ vertices. If $\asi(G) < \infty$, then
        \[
            \chi_\mathsf{B}(G) \,\leq\, 1 + \sum_{i = 1}^n \Delta(G[U_i]).
        \]
    \end{lemma}
    \begin{scproof}
        Let $G'$ be the spanning subgraph of $G$ with $E(G') \defeq E(G[U_1]) \sqcup \ldots \sqcup E(G[U_n])$. Applying Theorem~\ref{theo:indep_transversal} with $H = G_1 = G$ and $G_2 = G'$ yields a Borel $G$-independent complete section $S$ for $G'$. For each $1 \leq i \leq n$, let $U_i' \defeq U_i \setminus S$. Clearly, $\Delta(G[U_i']) \leq \Delta(G[U_i])$. Moreover, since $S$ meets every connected component of $G[U_i]$, 
         every component of $G[U_i']$ includes a vertex that is adjacent to $S$ in $G[U_i]$ and hence has degree at most $\Delta(G[U_i]) - 1$ in $G[U_i']$. Therefore, by a result of Conley, Marks, and Tucker-Drob \cite[p.~16, proof of Theorem 1.2]{MeasurableBrooks}, $\chi_\mathsf{B}(G[U_i']) \leq \Delta(G[U_i])$.

        To finish the argument, we apply Lemma~\ref{lemma:decomposition_simple} to the partition $V(G) = S \sqcup U_1' \sqcup \ldots \sqcup U_n'$ to get
        \[
            \chi_\mathsf{B}(G) \,\leq\, \chi_\mathsf{B}(G[S]) + \sum_{i=1}^n \chi_\mathsf{B}(G[U_i']) \,\leq\, 1 + \sum_{i=1}^n \Delta(G[U_i]). \qedhere
        \]
    \end{scproof}

    For the remainder of \S\ref{subsec:proof_Schreier}, we fix a group $\G$ generated by a finite symmetric set $F \subseteq \G \setminus \set{\mathbf{1}_\G}$. For ease of reference, let us restate Theorem~\ref{theo:BorelVizing} here:

    \begin{theocopy}{theo:BorelVizing}\label{theo:BorelVizing1}
        Suppose that $|F| \geq 2$ and for each $\sigma \in F$, the order of $\sigma$ is either even or at least $2^{30}\,|F|\log|F|$. Let $\bm{a} \colon \G \acts X$ be a free Borel action of $\G$ on a standard Borel space $X$ and let $G \defeq \mathsf{Sch}(\bm{a}, F)$ be the corresponding Schreier graph. If $\asi(G) < \infty$, then $\chi'_\mathsf{B}(G) \leq |F| + 1$.
    \end{theocopy}
    \begin{scproof}
        Recall that the line graph of $G$ is the graph $L(G)$ with vertex set $E(G)$ and edge set $\set{\set{e,e'} \in [E(G)]^2 \,:\, |e \cap e'| = 1}$. Note that we have $\chi'_\mathsf{B}(G) = \chi_\mathsf{B}(L(G))$. To bound $\chi_\mathsf{B}(L(G))$, we shall apply Lemmas~\ref{lemma:decomposition_simple} and \ref{lemma:decomposition_complicated} to certain partitions of $E(G)$.

    To begin with, we show that $\asi(L(G)) \leq \asi(G) < \infty$. Set $s \defeq \asi(G) < \infty$. By the Luzin--Novikov theorem (Theorem~\ref{theo:LN}), there is a Borel function $c \colon E(G) \to X$ such that $c(e) \in e$ for all $e \in E(G)$. Take any $R \in \N$ and let $X = U_0 \sqcup \ldots \sqcup U_s$ be a partition witnessing that $\si(G^{R+1}) \leq s$. Define a partition $E(G) = E_0 \sqcup \ldots \sqcup  E_s$ by $E_i \defeq c^{-1}(U_i)$. If $e$, $e' \in E(G)$ are joined by a path of length at most $R$ in $L(G)$, then $c(e)$ and $c(e')$ are joined by a path of length at most $R+1$ in $G$. It follows that the graphs $L(G)^R[E_i]$ are component-finite, and hence $\si(L(G)^R) \leq s$, as desired.
    
    Let $F_2$ be the set of all generators $\sigma \in F$ of order $2$ and let $F_{> 2} \subseteq F \setminus F_2$ be a subset formed by picking one member from each pair $\set{\sigma, \sigma^{-1}} \subseteq F \setminus F_2$. Then we have
    \[
        |F| \,=\, |F_2| + 2 |F_{> 2}|.
    \]
    We split the set $F_{>2}$ as $F_{>2} = F_{\mathsf{even}} \sqcup F_{\mathsf{odd}} \sqcup F_{\infty}$, where $F_\mathsf{even}$, $F_\mathsf{odd}$, and $F_\infty$ comprise the elements $\sigma \in F_{> 2}$ whose order is even, odd, or infinite respectively. For $\sigma \in F_2 \cup F_{> 2}$, define
    \[
        E_\sigma \,\defeq\, \set{\set{x, \sigma \cdot_{\bm{a}} x} \,:\, x \in X}.
    \]
    Note that the sets $(E_\sigma \,:\, \sigma \in F_2 \cup F_{> 2})$ form a partition of $E(G)$.
    
    Let $L_\sigma \defeq L(G)[E_\sigma]$. If $\sigma \in F_2$, then $E_\sigma$ is an $L(G)$-independent set. For $\sigma \in F_\mathsf{even} \cup F_\mathsf{odd}$, every component of $L_\sigma$ is a cycle of length $\mathrm{ord}(\sigma)$, and for $\sigma \in F_\infty$, every component of $L_\sigma$ is an infinite path. If $\sigma \not \in F_\infty$, then $L_\sigma$ is component-finite and hence, by Corollary \ref{corl:coloring}, $\chi_\mathsf{B}(L_\sigma) = \chi(L_\sigma)$. 

    Let $L_0$ and $L_1$ be the subgraphs of $L(G)$ induced by $\bigcup_{\sigma \in F_2 \cup F_\mathsf{even}} E_\sigma$ and $\bigcup_{\sigma \in F_\mathsf{odd} \cup F_\infty} E_\sigma$ respectively. By Lemma~\ref{lemma:decomposition_simple}, 
    \[
        \chi_\mathsf{B}(L_0) \,\leq\, \sum_{\sigma \in F_2 \cup F_\mathsf{even}} \chi_\mathsf{B}(L_\sigma) \,=\, |F_2| + 2|F_\mathsf{even}|.
    \]
    Note that $\Delta(L(G)) < 2|F|$. 
    Since every component of $L_\sigma$ for $\sigma \in F_\mathsf{odd} \cup F_\infty$ contains at least
    \[
        2^{30}\,|F|\log|F| \,\geq\, 2^{25} \cdot 2|F|\log(2|F|)
    \]
    vertices, we may apply Lemma~\ref{lemma:decomposition_complicated} to conclude that
    \[
        \chi_\mathsf{B}(L_1) \,\leq\, 1 + \sum_{\sigma\in F_\mathsf{odd} \cup F_\infty} \Delta(L_\sigma) \,=\, 1 + 2 |F_\mathsf{odd}| + 2|F_\mathsf{\infty}|.
    \]
    By Lemma~\ref{lemma:decomposition_simple} again,
    \[
        \chi_\mathsf{B}(L(G)) \,\leq\, \chi_\mathsf{B}(L_0) \,+\, \chi_\mathsf{B}(L_1) \,\leq\, |F_2| + 2|F_\mathsf{even}| + 1 + 2 |F_\mathsf{odd}| + 2|F_\mathsf{\infty}| \,=\, |F| + 1,
    \]
    and the proof is complete.
    \end{scproof}

    \subsection*{Acknowledgment}

    We are very grateful to the anonymous referee for their feedback.

	\printbibliography
    
\end{document}